\documentclass{article}

\usepackage{arxiv}

\usepackage[utf8]{inputenc}
\usepackage[T1]{fontenc}
\usepackage{hyperref}
\usepackage{url}
\usepackage{doi}

\usepackage[english]{babel}
\usepackage{amsmath}
\usepackage{amsfonts}
\usepackage{graphicx}
\usepackage{natbib}
\usepackage{enumitem}

\usepackage{amssymb}

\usepackage{color}
\usepackage{multirow}
\usepackage{mathtools}
\usepackage{pdflscape}
\usepackage{tabularx} 
\usepackage{mathrsfs} 
\usepackage{float} 
\usepackage{booktabs} 
\usepackage[algoruled,vlined,linesnumbered]{algorithm2e}
\usepackage{setspace}
\usepackage{threeparttable}
\usepackage{longtable}
\usepackage{subfigure}

\DeclarePairedDelimiter\abs{\lvert}{\rvert}

\allowdisplaybreaks

\newcommand{\cO}{\mathcal{O}}
\newcommand{\cM}{\mathcal{M}}
\newcommand{\cN}{\mathcal{N}}
\newcommand{\cF}{\mathcal{F}}
\newcommand{\cB}{\mathcal{B}}

\newcommand{\cP}{\mathcal{P}}
\newcommand{\cU}{\mathcal{U}}
\newcommand{\cA}{\mathcal{A}}

\newcommand{\MHa}{\texttt{ILS-Math\textsubscript{1}}}
\newcommand{\MHb}{\texttt{ILS-Math\textsubscript{2}}}
\newcommand{\MHc}{\texttt{ILS-Math\textsubscript{3}}}
\newcommand{\MHd}{\texttt{GRASP-Math\textsubscript{1}}}
\newcommand{\MHe}{\texttt{GRASP-Math\textsubscript{2}}}
\newcommand{\MHf}{\texttt{GRASP-Math\textsubscript{3}}}

\newcommand{\Rplus}{\protect\hspace{-.1em}\protect\raisebox{.35ex}{\smaller{\smaller\textbf{+}}}}
\newcommand{\Cpp}{\mbox{C\Rplus\Rplus}\xspace}

\newcommand{\ihat}{{\hat{\imath}}}

\allowdisplaybreaks

\title{Matheuristics for a Parallel Machine Scheduling Problem with Non-Anticipatory Family Setup Times: Application in the Offshore Oil and Gas Industry}

\author{
\href{https://orcid.org/0000-0003-3042-9249}{\includegraphics[scale=0.06]{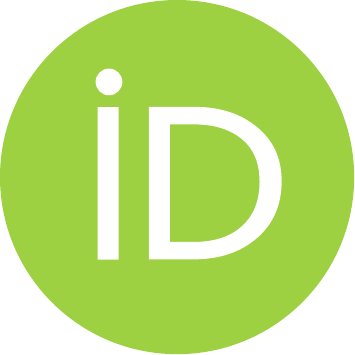}\hspace{1mm}Victor Abu-Marrul} \\
Departamento de Engenharia Industrial \\
Pontif\'{\i}cia Universidade Cat\'olica do Rio de Janeiro \\
\texttt{victorabu@aluno.puc-rio.br} \\
\And
\href{https://orcid.org/0000-0001-5715-149X}{\includegraphics[scale=0.06]{orcid.pdf}\hspace{1mm}Rafael Martinelli} \\
Departamento de Engenharia Industrial \\
Pontif\'{\i}cia Universidade Cat\'olica do Rio de Janeiro \\
\texttt{martinelli@puc-rio.br} \\
\And
\href{https://orcid.org/0000-0002-3980-0342}{\includegraphics[scale=0.06]{orcid.pdf}\hspace{1mm}Silvio Hamacher} \\
Departamento de Engenharia Industrial \\
Pontif\'{\i}cia Universidade Cat\'olica do Rio de Janeiro \\
\texttt{hamacher@puc-rio.br} \\
\And
\href{https://orcid.org/0000-0002-0078-3426}{\includegraphics[scale=0.06]{orcid.pdf}\hspace{1mm}Irina Gribkovskaia} \\
Faculty of Logistics, \\
Molde University College \\
\texttt{Irina.Gribkovskaia@himolde.no}
}

\renewcommand{\headeright}{A Preprint}
\renewcommand{\undertitle}{A Preprint}
\renewcommand{\shorttitle}{Matheuristics for a Parallel Machine Scheduling Problem}

\hypersetup{
pdftitle={Matheuristics for a Parallel Machine Scheduling Problem with Non-Anticipatory Family Setup Times: Application in the Offshore Oil and Gas Industry},
pdfsubject={math.OC},
pdfauthor={Victor Abu-Marrul, Rafael Martinelli, Silvio Hamacher, Irina Gribkovskaia},
pdfkeywords={Parallel machine scheduling, Family scheduling, Batch scheduling, Matheuristic, Offshore industry logistics, Ship scheduling},
}

\begin{document}
\maketitle

\begin{abstract}
In this paper, we address a variant of a batch scheduling problem with identical parallel machines and non-anticipatory family setup times to minimize the total weighted completion time. We developed an ILS and a GRASP matheuristics to solve the problem using a constructive heuristic and two MIP-based neighborhood searches, considering two batch scheduling mathematical formulations. The problem derives from a ship scheduling problem related to offshore oil \& gas logistics, the Pipe Laying Support Vessel Scheduling Problem (PLSVSP). The developed methods overcome the current solution approaches in the PLSVSP literature, according to experiments carried out on a benchmark of 72 instances, with different sizes and characteristics, in terms of computational time and solution quality. New best solutions are provided for all medium and large-sized instances, achieving a reduction of more than 10\% in the objective function of the best case.
\end{abstract}

\keywords{Parallel machine scheduling \and Family scheduling \and Batch scheduling \and Matheuristic \and Offshore industry logistics \and Ship scheduling}

\section{Introduction}
\label{sec:introduction}

The discovery of Brazilian pre-salt fields in 2006, the largest one in recent years, duplicated Brazilian oil and gas reserves. These fields are located in ultra-deep waters, below the ocean salt deposits, exceeding 2,000 meters of water depth. Exploration and Production (E\&P) are more challenging in this region in terms of technology and sustainability \citep{beltrao2009ss, allahyarzadeh2018energy, haddad2010economic}. Wells drilled in the Brazilian pre-salt basin are connected to surface platforms by flexible pipelines that better fits to the high water depths. Pipe Laying Support Vessels (PLSVs) are responsible for loading these pipelines at the port, transporting them to the wells site, laying them out in the ocean, and connecting them between the wells and platforms, allowing production to start \citep{speight2014handbook, clevelario2010flexible}. The \textit{PLSV Scheduling Problem} (PLSVSP) consists of servicing a demand of sub-sea oil wells connections, finding the best schedule for a limited PLSV fleet, prioritizing the completion of wells with higher production levels.

The PLSVSP has been addressed previously in the literature. \cite{queirozemendes2012} and \cite{bremenkamp2016} simplified the problem by grouping pipeline connections to compose jobs, modeling it as a classic identical parallel machine scheduling problem. In the recent works of \cite{cunhaILSmic}, \cite{cunhails}, \cite{AbuMArrul2020}, and \cite{Marrul2020}, the problem is formulated as an identical parallel machine scheduling problem with batching. \cite{cunhaILSmic} and \cite{cunhails} applied heuristics to reschedule a small set of instances with similar characteristics to minimize the impacts caused by disruptions in a given schedule. \cite{AbuMArrul2020} introduced three Mixed Integer Programming (MIP) formulations for the PLSVSP. They tested their approach on 72 instances generated from actual pre-salt data provided by a Brazilian company. A batch scheduling formulation using a dispatching rule to sequence operations within batches produced the best solutions, running for six hours. \cite{Marrul2020} developed several heuristics for building fast initial solutions for the PLSVSP. In our paper, we extend the works of \cite{AbuMArrul2020} and \cite{Marrul2020}, proposing a new variant on the batch scheduling formulation. We developed matheuristics using the constructive heuristic and the batch formulations to improve the solution quality and reduce computational time compared to previous works. An Iterated Local Search (ILS) and a Greedy Randomized Adaptive Search Procedure (GRASP) matheuristic for the PLSVSP are introduced, using two MIP-based neighborhood searches.

Matheuristics are hybrid approaches that combine concepts of metaheuristics and exact methods, being a growing field in the literature due to the improvement of computers and solvers \citep{thompson2018exact}. Some researchers applied matheuristics to solve machine scheduling problems. \cite{billaut2015single} handled a single machine environment using a two-step approach, combining a beam search algorithm with a MIP-based neighborhood search. 
Regarding parallel machines environments, \cite{ekici2019application} applied a Tabu Search matheuristic, prohibiting some job-machine assignments during its execution. \cite{ozer2019mip} developed a two-step approach combining a genetic algorithm with a MIP model. \cite{woo2018matheuristic} proposed a two-step approach, grouping jobs in so-called buckets using metaheuristics and assigning them to machines by a mathematical model. \cite{fanjul2017models} introduced matheuristics using constraint relaxation, limiting job-machine assignments, and using a MIP-based neighborhood search to optimize subsets of jobs. Other researchers applied matheuristics to solve flow-shop scheduling problems. \cite{ta2018matheuristic}, \cite{della2014matheuristic}, and \cite{della2019minimizing} used MIP-based neighborhoods in positional scheduling formulations, solving sub-problems for a limited number of positions. \cite{lin2016optimization} and \cite{lin2019makespan} converted the flow-shop problem into a traveling salesman problem, building an initial heuristic solution and solving a mathematical model. To our knowledge, \cite{monch2018matheuristic} is the only work that applies matheuristics to a batch scheduling formulation. However, they do not apply MIP-based neighborhood search but use a two-step approach that combines a genetic algorithm to compose batches and a mathematical model to assign and sequence them on the machines. Matheuristics have also been successfully applied to many scheduling problems with realistic backgrounds. See, for instance the works of \cite{martinelli2019strategic}, \cite{kalinowski2020scheduling}, and \cite{grenouilleau2020new}, related to the scheduling of mining activities, rail network maintenance, and home health care services, respectively.

The PLSVSP can be modeled as a variant of an identical parallel machine scheduling problem with non-anticipatory family setup times and batching. In this analogy, vessels are machines, jobs represent wells, and pipeline connections are the operations. The non-anticipatory family setup times represent pipeline loading times, while PLSV voyages from the ports to the wells site are the batches. In this machine environment, a set of parallel machines is available to perform a given set of tasks. Machines are called identical when tasks processing times are fixed and, therefore, machine-independent \citep{pinedo2012scheduling}. In most studies, tasks are called jobs, as each task requires only a single operation to complete. However, we represent tasks as operations, given that a set of operations must be finished to complete a job, in our problem. Besides, when family setup times are considered, tasks are grouped into families by similarity, and a setup time must be scheduled whenever a machine changes the execution of a task from one family to another. In this class of problems, the combination of one setup time and its subsequent tasks is called a batch. Also, when these setup times are non-anticipatory, the starting time of a batch cannot be anticipated, that is, it depends on the release of the tasks assigned to it. Parallel machine scheduling problems with family setup times are an important field in the scheduling literature. \cite{Shin2004}, \cite{Schaller2014}, \cite{Obeid2014}, and \cite{Ciavotta2016} applied heuristics to solve parallel machine scheduling problems with realistic backgrounds but not dealing with non-anticipatory setup times. The reader can be referred to \cite{Allahverdi2015} and \cite{Potts} for more details about scheduling problems with setup times considerations and batching. 

The novelty of our research relies on the development of two new MIP-based neighborhood searches for batch scheduling formulations, testing its efficiency in an ILS, and a GRASP matheuristic procedures. We also provide an extension of a batch scheduling formulation for the PLSVSP, considering the sequence of operations inside batches as a model's decision. Another contribution concerns the application to a problem with a realistic background, where the search for suitable solutions is crucial for sustainability in the company's operation. The method developed in this work outperforms the current PLSVSP solution approaches in terms of solution quality and computational time when solving a benchmark of 72 PLSVSP instances developed from real-life data. Besides, dealing with a complex problem allows the use of the presented concepts on other batch scheduling problems by adjusting some formulation constraints. In addition, we contribute to the increase in the literature on matheuristics applied to parallel machine scheduling problems with batching.

The remainder of this paper is organized as follows. A description of the problem and its correspondence as an identical parallel machine scheduling problem is given in Section \ref{sec:description}. Two mathematical formulations for the PLSVSP are provided in Section \ref{sec:formulation}. In Section \ref{sec:constructive}, a constructive heuristic for the PLSVSP is presented. The MIP-based neighborhood searches are introduced in Section \ref{sec:mipneig} with the ILS and GRASP matheuristics described in Section \ref{sec:matheuristics}. In Section \ref{sec:experiments}, computational experiments are presented and discussed. Finally, the conclusion and some perspectives for future research are provided in Section \ref{sec:conclusion}.

\section{Problem Description}
\label{sec:description}

The PLSVSP consists of scheduling a given PLSV fleet to meet a demand for pipeline connections in different sub-sea oil wells. The objective is to anticipate the completion of more productive wells (i.e., wells drilled in larger oil deposits). PLSVs are responsible for transporting the pipelines from the port to the wells site and connecting them between the wells and the production platforms. Each well includes a number of pipelines that must be connected to complete it, allowing its production to begin. Navigation times are disregarded in the problem due to the wells' proximity in the Brazilian pre-salt basin, reducing the importance of vessel routing, thus defining a scheduling problem. The PLSV planners must define schedules of voyages for the vessels in which each voyage consists of the loading process of the pipelines at the port, followed by one or more pipeline connections to perform. Three main decisions are made when scheduling a PLSV fleet:

\begin{enumerate}
    \item Definition of the voyages;
    \item Assignment of the voyages to the vessels;
    \item Sequencing of the voyages on each vessel.
\end{enumerate}

In the first step, the planners must consider the family of each operation to compose the voyages. Each family represents a group of operations with a similar loading process at the port, and a voyage must only include operations of the same family. Moreover, the loading duration of each voyage is family-dependent. On the assignment step, the schedulers must check if a vessel is eligible and has enough space on the deck to perform all pipeline connections defined for a voyage since the fleet is heterogeneous. Finally, on the voyage sequencing step, the planners must consider the arrival of the pipelines at the port, which may occur on different days, and the day on which each vessel is available for starting their activities. The loading process can only begin when all pipelines defined for a voyage are available at the port. Thus, the defined sequence directly affects the vessels' idleness and the delay in starting subsequent voyages. Due to the problem constraints and the concern with the solution's quality, all activities from this decision process are carried out simultaneously.

The PLSVSP is described in Section \ref{sec:pmach} as an identical parallel machine scheduling problem with non-anticipatory family setup times and batching. We provide a mapping between the machine scheduling aspects and the PLSVSP context in Appendix \ref{app:mapping} to help readers clarifying this relation.

\subsection{Identical Parallel Machine Scheduling Approach}
\label{sec:pmach}

The notation and assumptions considered in the PLSVSP are described below, according to the machine scheduling theory. 

\begin{enumerate}
  \item There is a set $\cO$ of operations, where each operation $i$ has a processing time $p_i$, a release date $r_i$, a load occupancy $l_i$, and a family $f_i$.
  \item There is a set $\cF$ of families where $s_f$ is the setup time for a family $f$. 
  \item All operations must be scheduled without preemption in a set $\cM$ of identical machines. 
   \item Each machine $k \in \cM$ is available to process operations from its release date $r_k$ and has a capacity $q_k$.  
  \item Machines are called identical due to the fixed processing times of the operations.
  \item A subset $\cM_i$ defines the eligible machines for executing each operation~$i$.~ Conversely, $\cO_k$ is a subset of operations $i \in \cO$ that a machine $k \in \cM$ is eligible to execute, i.e., $\cO_k = \{i \in \cO \: | \: k \in \cM_i\}$.
  \item A family setup time is incurred on three occasions: while changing the execution of operations from different families, before the first operation on each machine, or when the machine's capacity is reached.
  \item We define $Batch$ as a combination of one family setup time followed by a sequence with one or more operations from the same family.
  \item The batching mode is Serial-Batching. Thus, the processing time of a batch is given by the sum of the operations' processing times within the batch plus the setup time duration regarding the batch family.
  \item The size of a batch, computed by the sum of the load occupancy of the operations within it, must respect the capacity of the machine assigned to execute it. 
  \item The setup times are non-anticipatory, i.e., a $Batch$ can only start when all operations within it are released.
  \item There is a set $\cN$ of jobs where each job $j$ is associated with a subset $\cO_j$ of operations and has a weight $w_j$ defining its priority.
  \item We use $\cN_i$ to identify a subset of jobs associated with operation $i$, since in the PLSVSP, one operation might be related to several jobs simultaneously.
  \item A job is completed when all of its associated operations are concluded.  Thus, the completion time ($C_j$) of a job is the maximum completion time of the associated operations ($C_j = \max_{i \in \cO_j} C_i$). $C_i$ is the completion time of operation~$i$.
  \item The objective function is to minimize the total weighted completion time of all jobs, defined as $\sum_{j \in \cN} w_jC_j$.
  
\end{enumerate}

\begin{figure}[htbp] 
	\centerline{\includegraphics[scale=0.215, trim = 10 0 0 0, clip]{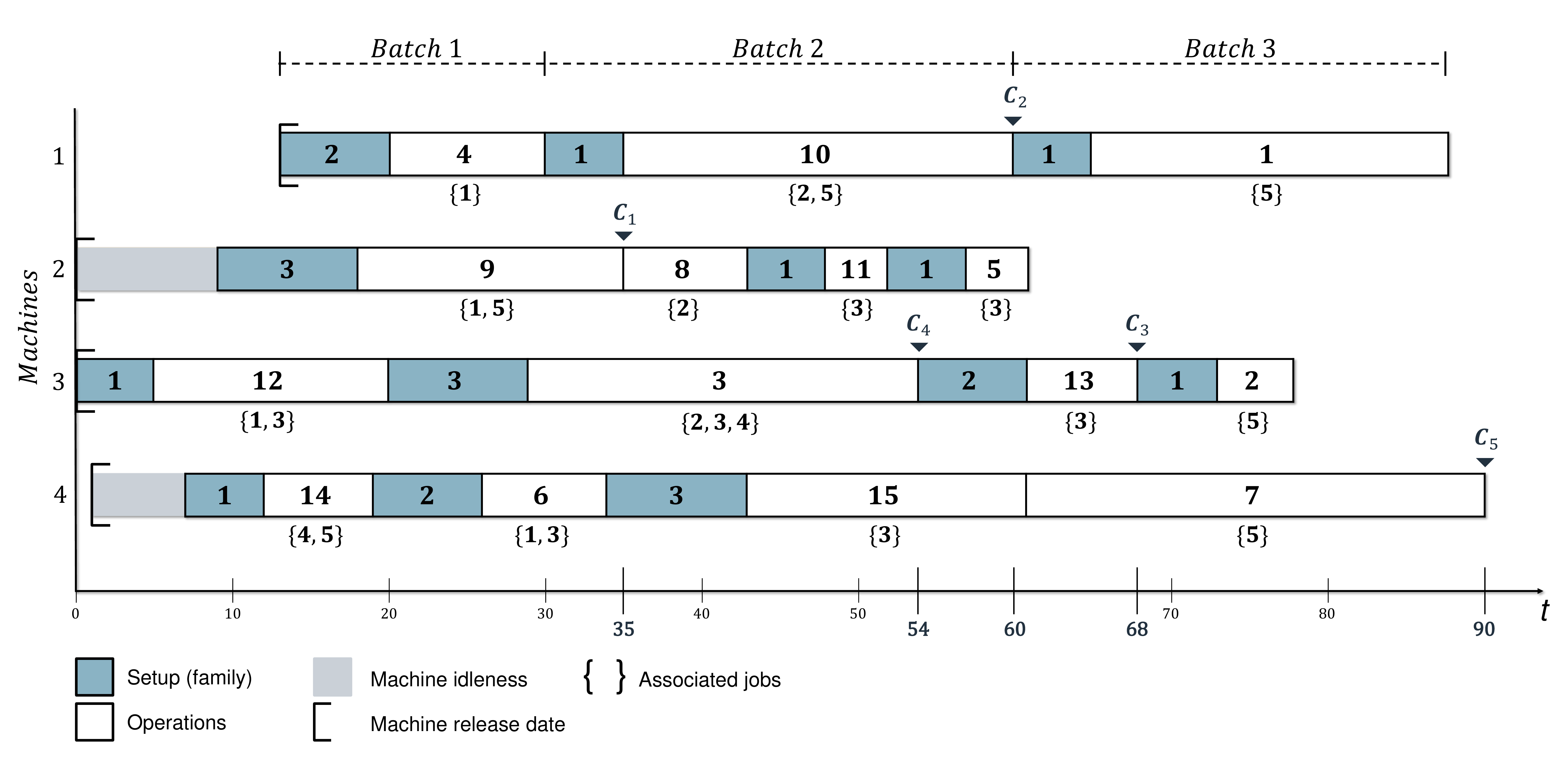}}
	\caption{PLSV scheduling example with 15 operations, 5 Jobs and 4 machines. Adapted from~\cite{AbuMArrul2020}. \label{fig:scheduling-example}}
\end{figure}

In Figure \ref{fig:scheduling-example}, a PLSV schedule example is depicted, consisting of 15 operations, assigned, sequenced, and forming batches on four machines. The batches on the first machine are depicted on the top. Setup times are depicted with its respective families. The jobs associated with each operation are shown below the allocations, and their completion times are highlighted on the horizontal axis (time horizon). The completion times of the jobs on the example are $C_1=35$, $C_2=60$, $C_3=68$, $C_4=54$ and $C_5=90$, with a total weighted completion time ($\sum_{j \in \cN} w_j C_j$) of $(35 \times 46) + (60 \times 40) + (68 \times 39) + (54 \times 13) + (90 \times 3) = 7,634$. The complete description of the example data is presented in Appendix \ref{app:example}.

\section{Mathematical Formulations}
\label{sec:formulation}

In this section, we present two MIP formulations for the PLSVSP. The first one is a batch formulation with a Weighted Short Processing Time (WSPT) dispatching rule \citep{pinedo2012scheduling}, introduced by \cite{AbuMArrul2020}. This formulation, which we refer to as \texttt{Batch-WSPT}, uses the WSPT dispatching rule to sequence operations heuristically inside the batches, reducing the solution space of the problem. To consider the entire solution space, we extend the \texttt{Batch-WSPT} formulation by adding inner batch sequencing variables. In this formulation, which we refer to as \texttt{Batch-S}, the sequence of operations inside batches is an optimization decision.

\subsection{Batch Formulation with WSPT Dispatching Rule}
\label{sec:formulation1}

To use the WSPT dispatching rule within the batch formulation, \cite{AbuMArrul2020} created the inner-batch precedence subsets $\cO_i$ (Equation~\ref{eq:heur-order}). Each subset $\cO_i$ is composed by the operations $\ihat \in \cO$ that will precede operation $i$ if both are scheduled in the same batch. These subsets consider an estimated weight $w_i$ for the operations to compute the WSPT rule ($w_i / p_i)$ and other breaking tie rules to achieve a complete ordering between operations. The complete ordering ensures that for any subset of operations assigned to the same batch, the model will be able to identify the sequence between them. To estimate the weight $w_i$ for each operation $i \in \cO$, we consider a proportional rule by splitting the weight of each job into equal parts for each of its associated operations, computed as $w_i=\sum_{j \in \cN_i} w_j / |\cO_j|$.

{\small
\begin{flalign}
\cO_i = & \Bigg\{ \ihat \in \cO \: \bigg| \: \left(                        \frac{w_\ihat}{p_\ihat} > \frac{w_i}{p_i}\right)    \vee \left(\frac{w_\ihat}{p_\ihat} = \frac{w_i}{p_i} \wedge w_\ihat > w_i\right) \vee 
        \left(\frac{w_\ihat}{p_\ihat} = \frac{w_i}{p_i} \wedge w_\ihat = w_i \wedge \ihat < i\right) \Bigg\} \:\:\:\: \forall i \in \cO &&
\label{eq:heur-order}
\end{flalign}
}
A set $\cB$ of batches is used in the formulation, and the subsets $\cB_k \subseteq \cB$ are generated to limit the number of available batches on each machine $k$. Taking into account the eligibility constraints, we limit $\abs{\cB_k} = \abs{\cO_k}$. Each element $b \in \cB_k$ defines the position of a batch in the schedule of machine $k$. When solving the PLSVSP, the model decides which batches to use, the family of these batches, and which operations to assign to it.

\bigskip
{\small
\noindent
The following binary variables are considered in the formulation:
\begin{align*}
X^{b}_{ik} =
& \begin{cases}
     1 & \text{if operation $i$ is scheduled in the $b$-th batch of machine $k$;} \\
     0 & \text{otherwise}. \\
   \end{cases} \\
Y^{b}_{fk} =
& \begin{cases}
     1 & \text{if the $b$-th batch of machine $k$ is of family $f$}; \\
     0 & \text{otherwise}. \\
   \end{cases}   
   &&&&&&&&&&&&&&&&&&&
\end{align*}
}

\noindent
Other continuous variables are considered in the formulation:
\medskip

{\small
\begin{itemize}
  \item{$S^b_k$ -- Starting time of the $b$-th batch on machine $k$.}
  \item{$P^b_k$ -- Processing time of the $b$-th batch on machine $k$.}
  \item{$C_i$ -- Completion time of operation $i$.}
  \item{$C_j$ -- Completion time of job $j$.}
\end{itemize}
}

\bigskip

\noindent
The \texttt{Batch-WSPT} formulation is as follows:

{\small
\begin{equation}
\min \sum\limits_{j \in \cN} w_j C_j 
\label{eq:sbf-obj}
\end{equation}
}
subject to

{\small
\begin{flalign}
\label{eq:sbf-1}
& \sum\limits_{k \in \cM_i}\sum\limits_{b \in \cB_k} X^{b}_{ik} = 1 & & \forall i \in \cO \\
\label{eq:sbf-2}
& \sum\limits_{f \in \cF} Y^{b}_{fk} \leq 1 & & \forall k \in \cM, b \in \cB_k \\
\label{eq:sbf-3}
& X^{b}_{ik} \leq Y^{b}_{f_i k} & & \forall i \in \cO, k \in \cM_i, b \in \cB_k\\
\label{eq:sbf-4}
& \sum\limits_{i \in \cO} l_i X^{b}_{ik} \leq q_k & & \forall k \in \cM, b \in \cB_k \\
\label{eq:sbf-8}
& P^b_k \geq \sum\limits_{i \in \cO} p_i X^{b}_{ik} + \sum\limits_{f \in \cF} s_f Y^{b}_{fk} & & \forall k \in \cM, b \in \cB_k \\
\label{eq:sbf-5}
& S^b_k \geq r_k & & \forall k \in \cM, b \in \cB_k \\
\label{eq:sbf-6}
& S^{b+1}_k \geq S^b_k + P^b_k & & \forall k \in \cM, b \in \cB_k \\
\label{eq:sbf-7}
& S^b_k \geq r_i X^{b}_{ik} & & \forall i \in \cO, k \in \cM_i, b \in \cB_k \\
\label{eq:sbf-9}
& C_i \geq S^b_k + p_i + s_{f_i} + \sum\limits_{\ihat \in \cO_i} p_{\ihat} X^{b}_{\ihat k} - \left(1 - X^{b}_{ik}\right)M & & \forall i \in \cO, k \in \cM_i, b \in \cB_k \\
\label{eq:sbf-10}
& C_j \geq C_i & & \forall j \in \cN, i \in \cO_j \\
\label{eq:sbf-13}
& C_i \geq 0 & & \forall i \in \cO \\
\label{eq:sbf-14}
& S^b_k, P^b_k \geq 0 & & \forall k \in \cM, b \in \cB_k\\
\label{eq:sbf-11}
& X^{b}_{ik} \in \{0, 1\} & & \forall i \in \cO, k \in \cM_i, b \in \cB_k \\
\label{eq:sbf-12}
& Y^{b}_{fk} \in \{0, 1\} & & \forall k \in \cM, b \in \cB_k, f \in \cF
\end{flalign}
}

The Objective Function~\eqref{eq:sbf-obj} minimizes the total weighted completion time of the jobs. Constraints~\eqref{eq:sbf-1} guarantees the execution of all operations. Constraints~\eqref{eq:sbf-2} limits a batch to be defined only for one family. Constraints~\eqref{eq:sbf-3} ensure that each operation follows the family defined for the assigned batch (family constraints). Constraints~\eqref{eq:sbf-4} limit the number of operations assigned to a batch considering the machine capacity (capacity constraints). Constraints~\eqref{eq:sbf-8} compute the processing time of batches. Constraints~\eqref{eq:sbf-5} limit the batch starting time according to the machine's release date. Constraints~\eqref{eq:sbf-6} ensure the batch starting time to respect the end of a previous batch. Constraints~\eqref{eq:sbf-7} force the batch starting time to respect the maximum release date between the operations scheduled in it (i.e., they guarantee the non-anticipatory setup time consideration). Constraints~\eqref{eq:sbf-9} compute the completion time of the operations. Constraints~\eqref{eq:sbf-10} compute the completion time of the jobs. Finally, Constraints~\eqref{eq:sbf-13}-\eqref{eq:sbf-12} present the variables' domains. Note that the domain of variable $X_{ik}^{b}$ guarantees the eligibility constraints by using the subsets $\cM_i$.

\subsection{Batch Formulation with Sequencing Variables}
\label{sec:formulation2}

As stated by \cite{AbuMArrul2020}, the \texttt{Batch-WSPT} does not consider the complete solution space of the problem since the sequence of operations inside batches is heuristically defined. To overcome this, we propose a variation on the \texttt{Batch-WSPT} formulation by adding a new variable and constraints to define the sequence of operations inside batches without considering the subsets $\cO_i$ (Equation~\ref{eq:heur-order}), named \texttt{Batch-S}.

To control the constraints and variables generation in the model, we consider the subsets $\mathscr{M}_{i \ihat}$ of machines that are eligible and with enough capacity for executing each pair ($i,\ihat \in \cO$) of operations in the same batch, defined as $\mathscr{M}_{i \ihat} = \{k \in (\cM_i \cap \cM_{\ihat}) \: | \: i \neq \ihat,\, f_i = f_{\ihat},\, l_i + l_{\ihat} \leq q_k \}$. A parameter $\mu_{i \ihat} \in \{0, 1\}$ is used to identify pairs of operations with at least one machine eligible to execute both in the same batch, equals 1 if $\abs{\mathscr{M}_{i \ihat}} > 0$, and zero otherwise. We do the same for operations triplets ($i, \ihat, i' \in \cO$). Thus, $\mathscr{M}_{i \ihat i'} = \{k \in (\cM_i \cap \cM_{\ihat} \cap \cM_{i'}) \: | \: i \neq \ihat, \ihat \neq i', i' \neq i,\, f_i = f_{\ihat} = f_{i'},\, l_i + l_{\ihat} + l_{i'} \leq q_k \}$, and $\mu_{i \ihat i'} \in \{0, 1\}$, equals 1 when $\abs{\mathscr{M}_{i \ihat i'}} > 0$, and zero otherwise.

\bigskip
\noindent
The following binary variable $Z_{i \ihat}$ is added to sequence operations inside batches:
{\small
\begin{align*}
Z_{i\ihat} =
& \begin{cases}
     1 & \text{if operation $i$ and $\ihat$ are scheduled in the same batch and $i$ precedes $\ihat$}; \\
     0 & \text{otherwise}. \\
   \end{cases}
   &&&&&&&&&&&&&&&&&&&&
\end{align*}
}

\noindent
The \texttt{Batch-S} formulation is as follows:

\begin{equation*}
\min \eqref{eq:sbf-obj}
\label{eq:min2}
\end{equation*}
subject to
\bigskip

\noindent
\eqref{eq:sbf-1}-\eqref{eq:sbf-7}, \eqref{eq:sbf-10}-\eqref{eq:sbf-14}
\vspace{-3.0mm}

{\small
\begin{flalign}
\label{eq:seq1}
& Z_{i\ihat} + Z_{\ihat i} \geq X^{b}_{ik} + X^{b}_{\ihat k} - 1 &  & \forall i, \ihat \in \cO, k \in \mathscr{M}_{i \ihat}, b \in \cB_k \: \\ 
\label{eq:seq2}
& Z_{i\ihat} + Z_{\ihat i} \leq 1 & & \forall i, \ihat \in \cO \: \big| \: \mu_{i \ihat} \\
\label{eq:seq3}
& Z_{i\ihat} + Z_{\ihat i'} + Z_{i' i} \leq 2 &  & \forall i, \ihat, i' \in \cO \: \big| \: \mu_{i \ihat i'} \\
\label{eq:seq4}
& C_i \geq S^b_k + p_i + s_{f_i} + \sum\limits_{\ihat \in \cO} p_{\ihat} Z_{\ihat i} - \left(1 - X^{b}_{ik}\right)M & & \forall i \in \cO, k \in \cM_i, b \in \cB_k \\
\label{eq:seq5}
& Z_{i \ihat} \in \{0, 1\} & & \forall i, \ihat \in \cO \: \big| \: \mu_{i \ihat}
\end{flalign}
}
Constraints \eqref{eq:seq1} identify when operations $i$ and $\ihat$ are scheduled in the same batch. Constraints \eqref{eq:seq2} ensure that only one of the variables that define the precedence between operations $i$ and $\ihat$ inside a batch will be considered. Constraints \eqref{eq:seq3} guarantee a complete ordering between operations inside a batch. Constraints \eqref{eq:seq4} replace Constraints \eqref{eq:sbf-9}. The completion time of an operation is now computed with the sequencing variable $Z_{i\ihat}$. Finally, Constraints \eqref{eq:seq5} present variable $Z_{i \ihat}$ domains.

\section{Constructive Heuristic}
\label{sec:constructive}

We use a constructive heuristic, introduced by \cite{Marrul2020}, to build the initial solutions for the PLSVSP. The authors tested 19 heuristics, combining several dispatching rules and ways of estimating weights for the operations, showing that one of them outperformed the others. The heuristic named \texttt{WMCT-WAVGA}, is an algorithm that, at each iteration, chooses an operation according to a Weighted Minimum Completion Time (WMCT) dispatching rule, scheduling it in the machine that minimizes the total weighted completion time in the partial schedule. Let $C_k$ be the completion time of each machine $k$, and $T_i$ be the minimum completion time among the set $\cM_i$ of eligible machines for each operation $i$, computed as $T_i = min_{k \in \cM_i}C_k$. Then, at each iteration, the operation with the highest priority value $\pi_i$ is selected, where $\pi_i$ is computed according to Equation \eqref{wmct}. The weight $w_i$ of operation $i$ is estimated as $w_i=\sum_{j \in \cN_i} w_j / |\cU_j|$, where $\cU$ is the set of unscheduled operations, and $\cU_j \subseteq \cU$ is a subsets of unscheduled operations associated with job $j$. In the first iteration, $\cU$ is equal to the set $\cO$ of operations. During the heuristic execution, the weights of unscheduled operations are adjusted at each iteration.

{\small
\begin{flalign}
\label{wmct}
& \pi_i = \frac{w_i}{\max(T_i, r_i) + p_i + s_{f_i}} \:\:\:\: \forall i \in \cU
\end{flalign}
}

The method creates batches by assigning operations sequentially to the machines. Therefore, at each iteration, the algorithm selects the next operation to schedule and assigns a machine to the operation. The method then decides whether to insert the selected operation in the last batch (called current batch) or to create a new batch to insert it in the assigned machine. If a new batch is created, the selected operation is sequenced as the first one inside the new batch, i.e., after a new family setup time is also inserted in the machine schedule. Otherwise, the operation is scheduled as the last one in the current batch on the assigned machine. There are three situations that force the creation of a new batch on the selected machine: (1) when the machine is empty; (2) when the current batch on the machine is of a different family from the selected operation; (3) when the insertion of the operation in the current batch exceeds the machine capacity. The \texttt{WMCT-WAVGA} heuristic, presented in Algorithm \ref{alg.constructive}, defines a list of schedules $\sigma = (\sigma_1, \ldots, \sigma_k)$ containing operations and families for each machine $k$. The families represent the setup times, indicating the beginning of a new batch.

\renewcommand{\gets}{\leftarrow}
\begin{algorithm}[htb!]
	\footnotesize
	\SetAlgoLined
	
	$C_k \gets r_k$, $S_k \gets r_k$, $L_k \gets 0$, $F_k \gets 0$, $\cA_k \gets \emptyset$, $\sigma_{k} \gets \emptyset$ \:\: $\forall k \in \cM$;
	
	$C_i \gets \infty \:\: \forall i \in \cO$;
	
	$\mathcal{U} \gets \cO$;

	\While{$\cU \neq \varnothing$}{
		
		$w_i \gets \sum\limits_{j \in \cN_i}\frac{w_j}{|\cU_j|}$, $T_i \gets \min\limits_{k \in \cM_i}C_k$, $\pi_i = \frac{w_i}{\max(T_i, r_i) + p_i + s_{f_i}}$ \:\:$\forall i \in \mathcal{U}$;
		
		Select operation $i^* \in \mathcal{U}$ that maximizes $\pi_i$;

		$\Delta_{i^* k} \gets \max(0, r_{i^*} - S_k)  \:\: \forall k \in \cM_{i^*}$;

		$C^{CB}_{i^* k} \gets C_k + \Delta_{i^* k} + p_{i^*}  \:\: \forall k \in \cM_{i^*}$;

		$C^{NB}_{i^* k} \gets \max(r_{i^*}, C_k)+s_{f_{i^*}} + p_{i^*}  \:\: \forall k \in \cM_{i^*}$;

		
		$ \mathcal{CB} \gets \bigg\{ cb_{i^* k} = w_{i^*} C^{CB}_{i^* k} + \sum\limits_{\ihat \in \mathcal{B}_k }{w_\ihat} \Delta_{i^* k} \:\: \big| \:\: k \in \cM_{i^*}, \: F_k=f_{i^*}, \: L_k+l_{i^*} \leq q_k \bigg\}$;
		
		$ \mathcal{NB} \gets \bigg\{nb_{i^* k} =w_{i^*} C^{NB}_{i^* k}\:\: \big| \:\: k \in \cM_{i^*}\bigg\}$;
		
		$b_{min} \gets \min\{b :  b \in (\mathcal{CB} \cup \mathcal{NB})\}$;
		
		Select $k^*$ corresponding to $b_{min}$;
		
		\eIf{$b_{min} \in \mathcal{CB}$}{

			$S_{k^*} \gets \max(r_{i^*}, S_{k^*})$; $C_{k^*} \gets C_{i^* k^*}^{CB}$;
			
			$L_{k^*} \gets L_{k^*} + l_{i^*}$; $\cA_{k^*} \gets \cA_{k^*} \cup \{i^*\}$;

		}{    
			
			$S_{k^*} \gets \max(r_{i^*}, C_{k^*})$; $C_{k^*} \gets C_{i^* k^*}^{NB}$;
			
			$L_{k^*} \gets l_{i^*}$; $\cA_{k^*} \gets \{i^*\}$;
			
			$\sigma_{k^*} \gets  \sigma_{k^*} \cup \{f_{i^*}\}$;
			
		}
		
		$\sigma_{k^*} \gets  \sigma_{k^*} \cup \{i^*\}$; $F_{k^*} \gets f_{i^*}$; $C_{i^*} \gets C_{k^*}$; $\mathcal{U} \gets \mathcal{U} \setminus \{i^*\}$;
		
	}
	
	\textbf{return} $\sigma$;
	
	\caption{\texttt{WMCT-WAVGA}} 
	\label{alg.constructive}
\end{algorithm}

New variables and sets are used to store information regarding the solution during the heuristic procedure. We use $S_k$, $L_k$, $F_k$, and $\cA_k$ to save the starting time, cumulative load, family, and set of assigned operations regarding the current batch on machine $k$. $\Delta_{ik}$ computes the delay of the starting time for the current batch on machine $k$ due to the insertion of operation~$i$ \textit{(caused by the operation's release date and the non-anticipatory setup time consideration)}. $C_{ik}^{CB}$ and $C_{ik}^{NB}$ represent the completion time of operation $i$, if inserted in the current batch or in a new batch on machine $k$, respectively. Finally, we consider two sets to store the feasible assignments and their costs: (1)~$\mathcal{CB}$  stores the set of feasible assignments $cb_{ik}$ of operation $i$ into the current batch on machine $k$ \textit{(to be feasible, the selected operation must respect the eligibility, family, and capacity constraints)}; (2)~$\mathcal{NB}$ stores the set of feasible assignments $nb_{ik}$ of operation $i$ into a new batch on machine $k$ \textit{(to be feasible, the selected operation must respect the eligibility constraint)}.

The algorithm starts by initializing the variables and sets (Lines 1-3). The method runs until all operations are scheduled (Line 4). The main loop (Lines 4-23) starts by computing the priority value ($\pi_i$) for all unscheduled operations (Line 5). In Line 6, the next operation~$i^*$ to schedule is selected according to the value of $\pi_i$. The delay at the starting time of the current batch on machine~$k$ is computed in Line 7. In Lines 8 and 9, the selected operation's completion times are computed, considering its insertion in the current batch or in a new batch on the assigned machine, respectively. Next, the method creates sets of feasible assignments (Lines 10 and 11). In Line 12, the assignment with minimum cost is selected, identifying the best machine~$k^*$ in Line 13. If the selected element belongs to the set $\mathcal{CB}$, operation~$i^*$ is inserted in the current batch on machine~$k^*$ (variables and sets are updated in Lines 15, and 16). Otherwise, a new batch on machine~$k$ is created for operation~$i$ (variables and sets are updated in Lines 18-20). The remaining variables and sets of the algorithm are updated in Line 22. Finally, the procedure returns the final solution $\sigma$ in Line 25.

Note that there are no rules to control the sequence of operations inside the batches during the \texttt{WMCT-WAVGA} heuristic execution. Since the \texttt{Batch-WSPT} formulation (Section~\ref{sec:formulation1}) generates these sequences heuristically, we need to update the subsets $\cO_i$ to use the constructed solutions in this formulation. Thus, we introduce a variable $\vartheta{i \ihat} \in \{0, 1\}$ which is equals 1 if a pair of operations ($i, \ihat \in \cO$) is scheduled in the same batch and $i$ precedes $\ihat$ in the constructed solution, and zero otherwise. Thus, the new definition of the subsets $\cO_i$ is shown in Equation~\eqref{eq:heur-order2}. Suppose that the solution depicted in Figure \ref{fig:scheduling-example} was generated by the \texttt{WMCT-WAVGA} heuristic. Then, the schedule $\sigma_k$ for the fourth machine ($k = 4$) would be $\sigma_{4} = \{f_{14}, 14, f_6, 6, f_{15}, 15, 7\}$, where $f_{14}$, $f_{6}$, and $f_{15}$ indicates families 1, 2, and 3, respectively. In that case, $\vartheta{\ihat i} = 1$ only for $\ihat = 15$ and $i = 7$, since this pair of operations are scheduled together in the last batch on machine~$4$.

{\small
\begin{flalign}
\cO_i = & \Bigg\{ \ihat \in \cO \: \bigg| \: \vartheta{\ihat i} \vee \left(\vartheta{i \ihat} = \vartheta{\ihat i} \wedge                        \frac{w_\ihat}{p_\ihat} > \frac{w_i}{p_i}\right)   
      \vee \left(\vartheta{i \ihat} = \vartheta{\ihat i} \wedge \frac{w_\ihat}{p_\ihat} = \frac{w_i}{p_i} \wedge w_\ihat > w_i\right) & & \nonumber \\
      & \vee 
        \left(\vartheta{i \ihat} = \vartheta{\ihat i} \wedge \frac{w_\ihat}{p_\ihat} = \frac{w_i}{p_i} \wedge w_\ihat = w_i \wedge \ihat < i\right) \Bigg\} \:\:\:\: \forall i \in \cO & &
\label{eq:heur-order2}
\end{flalign}
}

In the next section, we present the MIP-based neighborhood searches to consider in the local search step of our methods.

\section{MIP-based Neighborhood Searches}
\label{sec:mipneig}

We consider two MIP-based neighborhood searches, named \textit{Batch Windows} and \textit{Multi-Batches Relocate}, making use of the batch formulations presented in Section~\ref{sec:formulation}, to decompose the PLSVSP into smaller problems that can be optimized more quickly than the complete problem. The idea is to limit the number of integer variables to optimize at each iteration, fixing the remaining ones from a feasible solution. As mentioned before, we consider the initial feasible solution, the one provided by the \texttt{WMCT-WAVGA} heuristic (Algorithm~\ref{alg.constructive}). The approaches are detailed in the next sections.

\begin{figure}[htpb] 
	\centerline{\includegraphics[scale=0.175, trim = 0 0 0 10, clip]{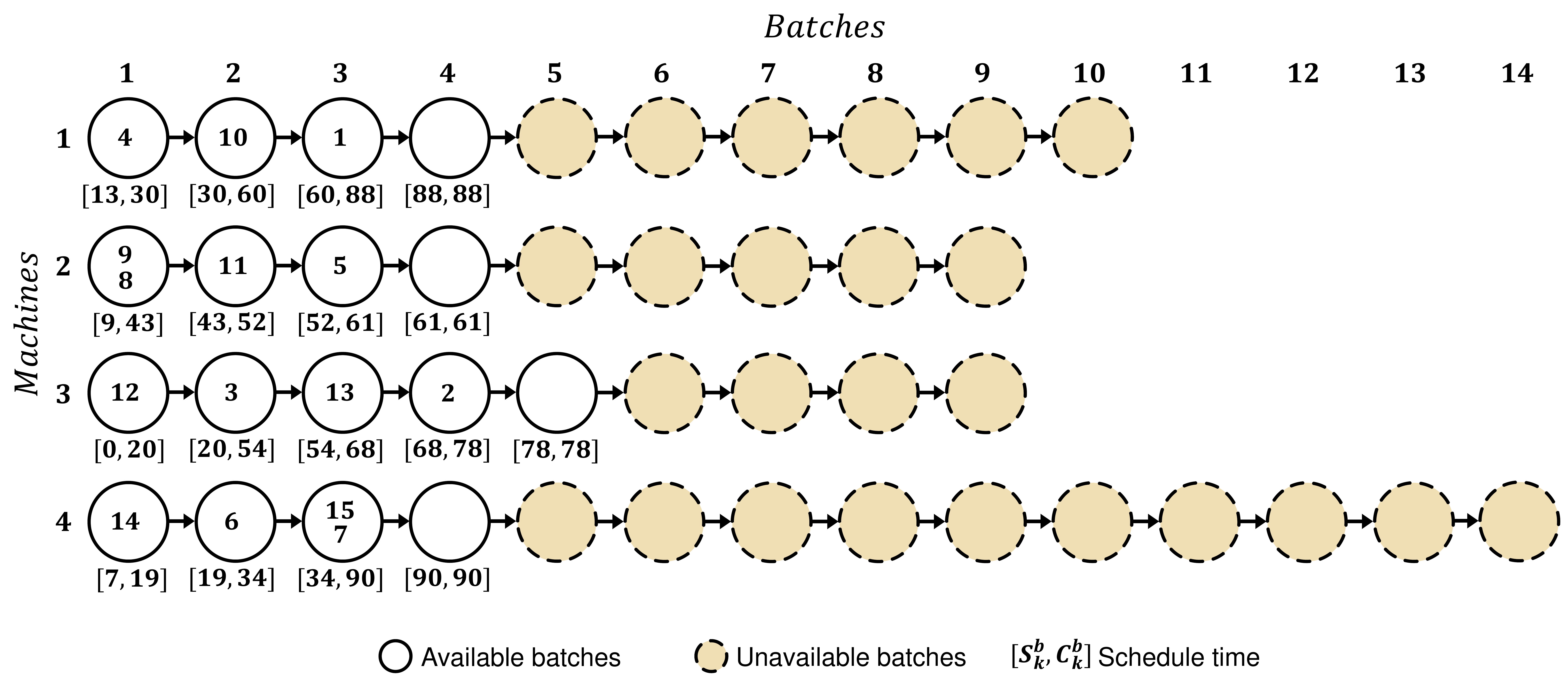}}
	\caption{Graph representation with batches of the PLSV scheduling example. \label{fig:graph}}
\end{figure}

In both methods, we use variables ${MB}_k \in \mathbb Z_{+}$ (Machine Batches) to limit the subset of batches to be considered on each machine $k \in \cM$ during the search. These variables are bounded by the size of the subsets $\cB_k$ on each machine $k$, thus $0 < {MB}_k \leq \abs{\cB_k}$. Throughout the search, only batches with position $b \leq {MB}_k$ are available to optimize. The remaining batches are fixed without any operations inside them. To help to clarify the idea, we show in Figure~\ref{fig:graph} a graph representation with batches of the PLSV schedule illustrated in Figure~\ref{fig:scheduling-example} (Section \ref{sec:description}). Each node represents a batch $b$ on a machine $k$, with the respective operations assigned to it described inside. Below each node, we show the starting time ($S_k^b)$ and the completion time ($C_k^b)$ of the respective batch. Note that the maximum number of batches on each machine $k$ is given by the total number of eligible operations. For instance, in machine~$1$, we have $|\cB_k|=10$ (batches 1 to 10), following the subsets $\cO_k$ defined in Table \ref{tab:machparameters} (Appendix \ref{app:example}), with $\cO_1=\{1,2,4,5,7,9,10,11,12,14\}$. Nodes are labeled as: (1)~Available batches (batches that can be used by the search procedures); (2)~Unavailable batches (batches generated by the formulation but not available for the search procedure in a given iteration, according to variables ${MB}_k$). The variables ${MB}_k$ are updated during the procedures, allowing an unavailable batch to become available at some point in the search.

The $MB_k$ variables are initialized according to a solution provided by the \texttt{WMCT-WAVGA} heuristic (Algorithm \ref{alg.constructive}), considering the number of used batches plus one extra batch on each machine $k$. Let $\eta_k^b \in \{0, 1\}$ be a variable equal 1 if a batch $b$ on a machine $k$ is used in a given solution (i.e., contains at least one operation scheduled in it), and zero otherwise. Then, $MB_k$ are initialized according to Equation~\eqref{eq:mbcalc2}. Note that, on machine 1, we have four available batches (1 to 4), with ${MB}_1=4$, although the solution, depicted in Figures \ref{fig:scheduling-example} and \ref{fig:graph}, only uses three batches on this machine. At each iteration of the MIP-based neighborhood searches, we check whether all available batches are being used on each machine $k$, updating $MB_k$ with Equation~\eqref{eq:mbcalc2}, if $\sum_{b \in \cB_k}{\eta^{b}_{k}}= MB_k$.

{\small
\begin{align}
\label{eq:mbcalc2}
& MB_{k} = 1 + \sum_{b \in \cB_k}\eta^{b}_{k} \:\:\:\: \forall k \in \cM
\end{align}
}

\subsection{Batch Windows}
\label{sec:batchwindows}

In this approach, we limit the subset of batches to optimize based on a defined time range interval. At each iteration, we only optimize batches that are scheduled inside the range. The Range Size ($RS$) is defined as a fraction of the makespan ($C_{max}~=~\max_{i \in \cO}C_i$) based on a given PLSV solution, computed as $RS = \lceil {\rho} \times C_{max} \rceil$, where $\rho \in [0, 1]$ is a parameter that defines the proportion of the makespan to consider. The complete search moves the optimization range from the end of the schedule to its beginning, with a step half the size of $RS$, ensuring overlap between iterations and that all batches are optimized at least once, totalizing in $\lceil C_{max}/(RS/2)\rceil$~-~1 iterations. We chose to move the search from the end of the schedule to its beginning based on preliminary experiments that showed advantages in this approach. We use $R_{begin}$ and $R_{end}$ to identify the beginning and the end of the optimization range to consider at each iteration.

\begin{figure}[http!]
	\centering
    \subfigure[Iteration 1: $R_{begin}=60$ and $R_{end}= 90$.\label{fig:batch_windows_1}]{%
		\includegraphics[scale=0.147, trim = 10 0 10 10, clip]{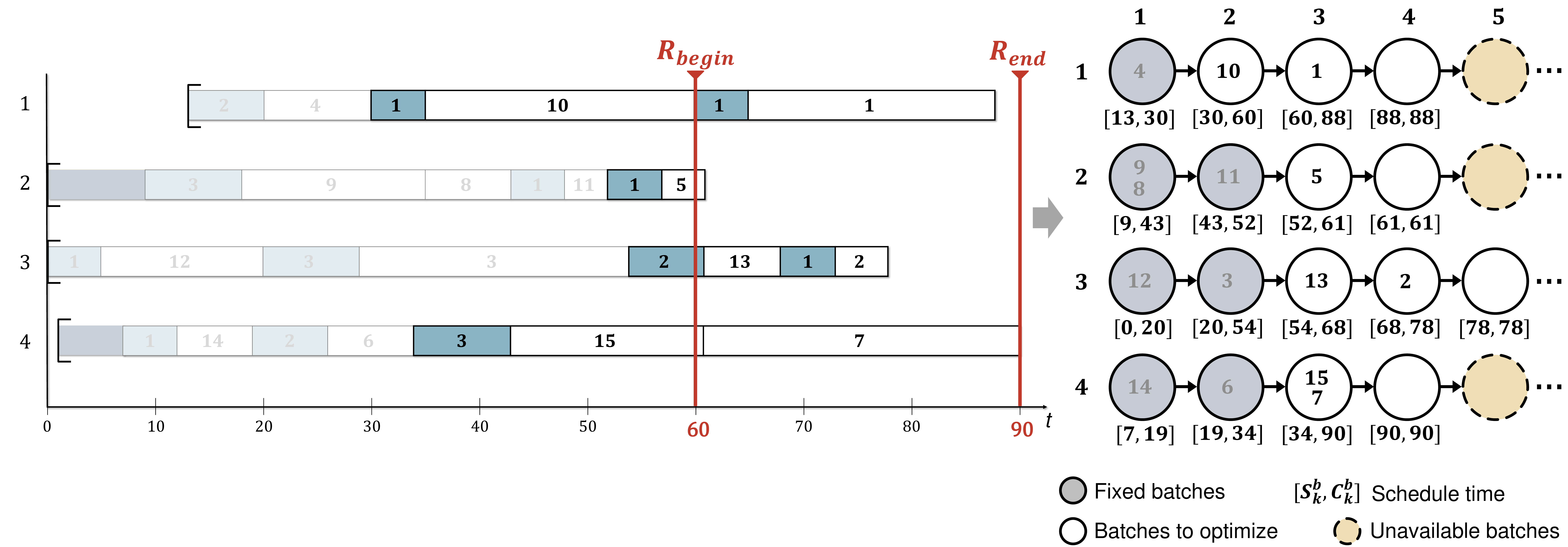}
	}
	\hfill
	\subfigure[Iteration 3: $R_{begin}=30$ and $R_{end}=60$.\label{fig:batch_windows_2}]{%
		\includegraphics[scale=0.147, trim = 10 0 10 10, clip]{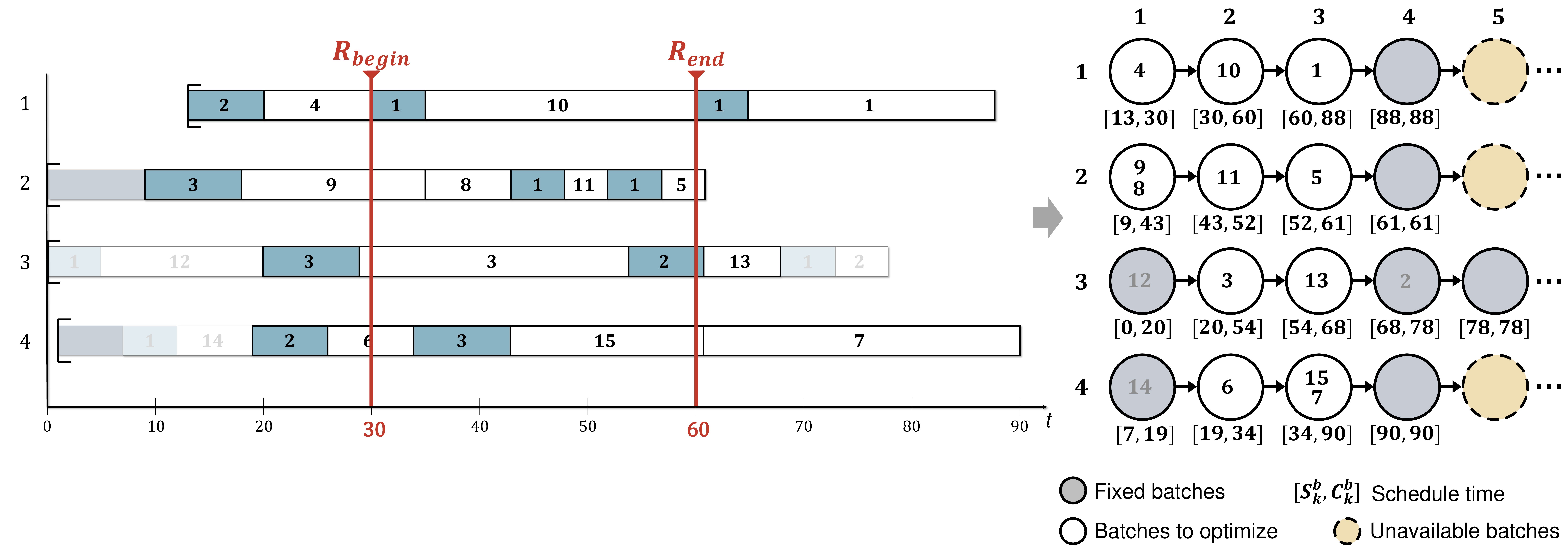}
	}
	\hfill
	\subfigure[Iteration 5: $R_{begin}=0$ and $R_{end}=30$.\label{fig:batch_windows_3}]{%
		\includegraphics[scale=0.147, trim = 10 0 10 10, clip]{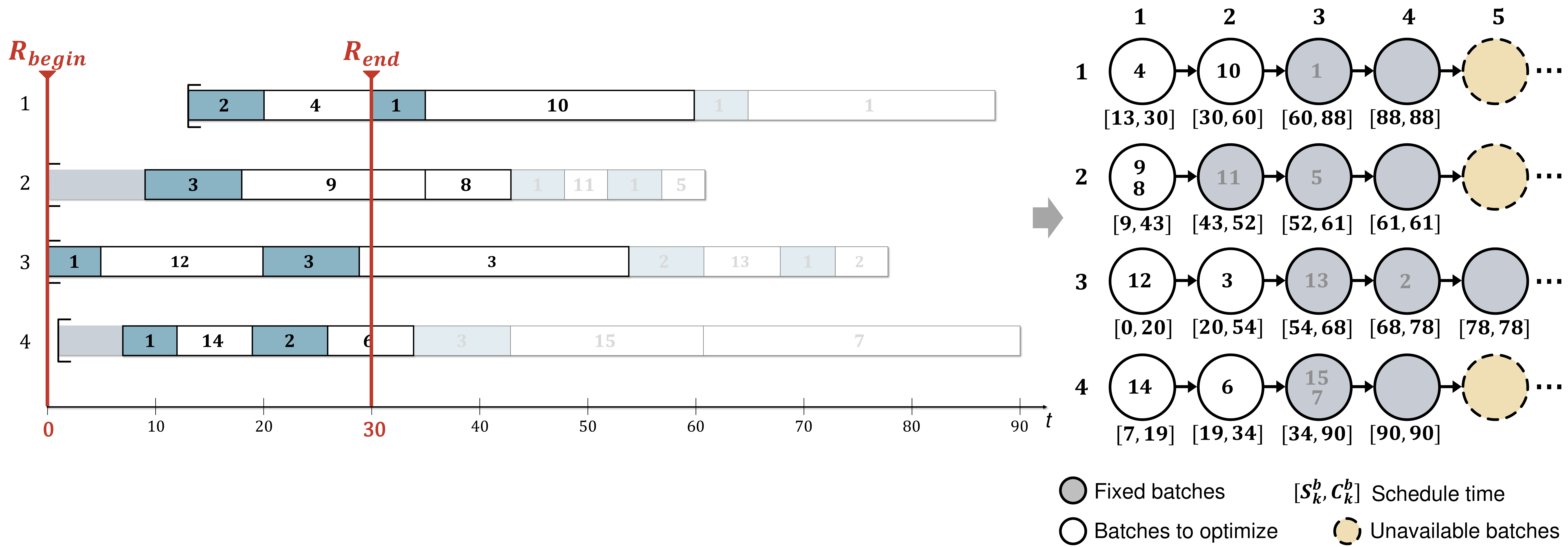}
	}
	\caption{Example of three iterations in the Batch Windows neighborhood search, showing a PLSV schedule and the optimization range on the left side, and the graph representation with the batches to optimize highlighted on the right side.}
	\label{fig:batchwindows}
\end{figure}

Given the PLSV schedule example shown in Figure~\ref{fig:scheduling-example} (Section~\ref{sec:description}), we depict in Figure~\ref{fig:batchwindows} how the optimization range defines the batch windows to optimize on each machine. We consider an optimization range of size 30 ($RS = 30)$, which would result in a total of five iterations, with the following ranges ($R_{begin}, R_{end}$): Iteration~1~(60,~90); Iteration~2~(45,~75); Iteration~3~(30,~60); Iteration~4~(15,~45); Iteration~5~(0,~30). To save space, we only show iterations 1, 3, and 5, depicted in Figures \ref{fig:batch_windows_1}, \ref{fig:batch_windows_2}, and \ref{fig:batch_windows_3}, respectively. In this example, we suppose the solution does not change during the search. On the right side of the figures, we show the solution graph representation with the batches on each machine, highlighting the ones to optimize at each iteration. Nodes are labeled as: (1)~Fixed batches (batches not selected in the given iteration); (2)~Batches to optimize (batches selected to optimize in the given iteration); (3)~Unavailable batches (batches generated by the formulation but not available for the search procedure in the given iteration, according to variables ${MB}_k$).

Note that a given optimization range defines different sizes of batch windows to optimize, on each machine, due to the continuous variables that compute the starting time ($S_k^b)$ and the completion time ($C_k^b)$ of each available batch. For instance, in Iteration 3, depicted in Figure \ref{fig:batch_windows_2}, on machines 1 and 2, there are three batches (1 to 3) to optimize, while on machines 3 and 4, there are two batches (2 and 3) to optimize. Let $\cO' \subseteq \cO$ be the subset of operations assigned to the selected batches to optimize in a given iteration. Then, $\cO'=\{1,3,4,5,6,7,8,9,10,11,13,15\}$ in iteration 3. The pseudo-code of the \textit{Batch Windows} is shown in Algorithm~\ref{alg.batchwindows}.

\renewcommand{\gets}{\leftarrow}
\begin{algorithm}[hbt!]
	\footnotesize
	\SetAlgoLined
	
	$C_{max} \gets \max_{i \in \cO}C_i$, where $C_i$ is given by the PLSVSP solution~$s$;
	
	$RS \gets \lceil {\rho} \times C_{max} \rceil$;
	
	$R_{begin} \gets \infty$;
	
	$R_{end} \gets C_{max}$;
	
	\While{$R_{begin} > 0$}{
		
		$C^b_k \gets S^b_k + P^b_k, \:\: \forall k \in \cM, b \in \cB_k$;
		
		$R_{begin} \gets \max(0, R_{end} - RS)$;
		
		Create the subset~$\cO'$ of operations assigned to batch $b$ on machine $k$ in which $S^b_k \leq R_{end}$ \textit{and} $C^b_k \geq R_{begin}$ \textit{and} $b \leq MB_k$ in solution~$s$;
		
		Solve \textit{Batch Formulation}, starting from solution $s$, for a subset of variables $X^{b}_{ik}$ in which $S^b_k \leq R_{end}$, $C^b_k \geq R_{begin}$, $b \leq MB_k$ and $i \in \cO'$;
		
		$R_{end} \gets R_{begin} + RS / 2$;
		
		Update variables $MB_k$ according to Equation \ref{eq:mbcalc2}, if $\sum_{b \in \cB_k}{\eta^{b}_{k}}= MB_k$;
		
	}
	
	\caption{Batch Windows ($s, \rho$, $MB_k$)} 
	\label{alg.batchwindows}
\end{algorithm}

Algorithm \ref{alg.batchwindows}  starts by computing values for $C_{max}$, $RS$, $R_{begin}$, and $R_{end}$ (Lines 1-4), according to a given solution $s$ and the defined parameter $\rho$. The main loop of the algorithm (Lines 5-12) is repeated until $R_{begin}$ reaches zero. The completion time ($C^b_k$) of each batch $b$ on each machine $k$ is computed in Line 6. The algorithm updates $C^b_k$ and $R_{begin}$ (Lines 6 and 7), defining the subset $\cO'$ in Line 8. The \textit{Batch Formulation} is solved for all variables $X^{b}_{ik}$ of batches scheduled within the optimization range (Line 9). After solving the model, the value of $R_{end}$ is updated for the next iteration in Line 10. Finally, the procedure to update the number of available batches on each machine is executed in Line 11.

\subsection{Multi-Batches Relocate}
\label{sec:multibatchrelocate}

In this approach, we randomly select the batches to optimize at each iteration, not allowing the selection of batches already optimized in previous iterations. The complete search ends when each batch is optimized exactly once in one of the iterations. We compute the Number of Batches ($NB$) to optimize at each iteration based on a given parameter $\varphi \in [0, 1]$, as $NB = \lceil {\varphi} \times \sum_{k \in \cM}{MB_k} \rceil$. The method runs with a total of $\lceil \sum_{k \in \cM}{MB_k} / NB \rceil$ iterations.

\begin{figure}[hbt!] 
	\centerline{\includegraphics[scale=0.143, trim = 0 0 0 0, clip]{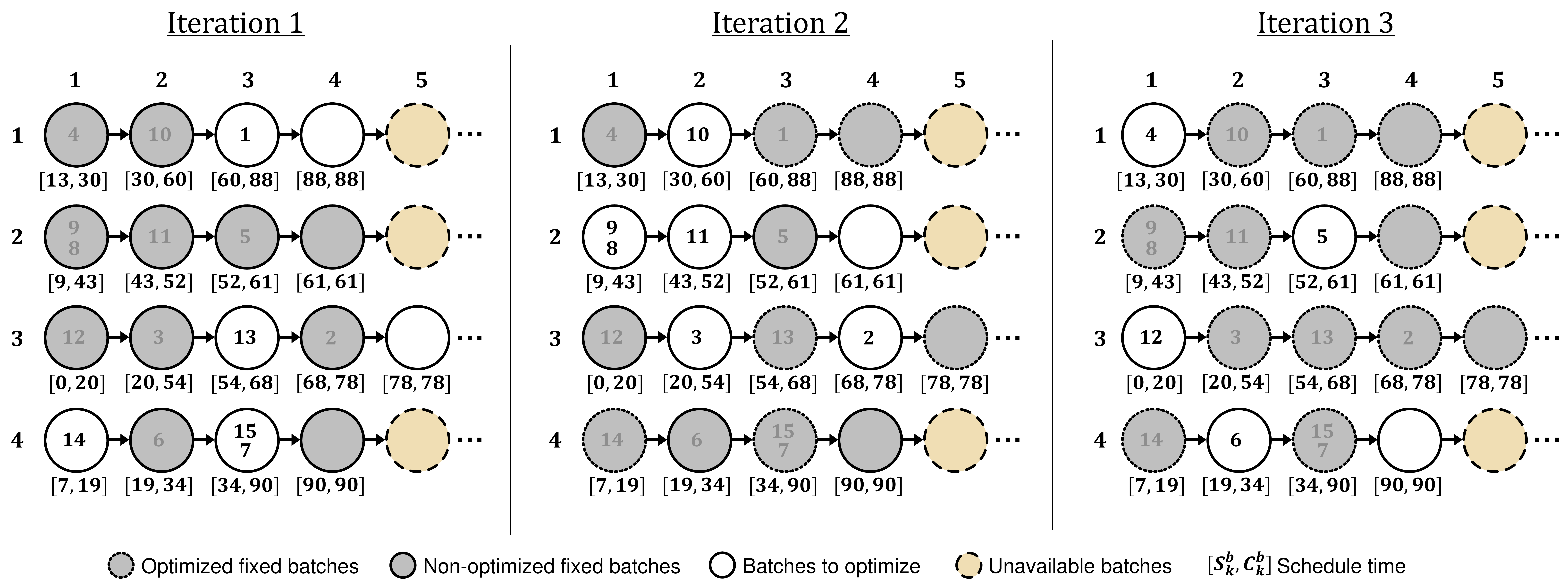}}
	\caption{Example with three iterations of batches selection in the Multi-Batches Relocate. 
	\label{fig:multirelocate}}
\end{figure}

An example of the \textit{Multi-Batches Relocate} is depicted in Figure~\ref{fig:multirelocate}, using the graph representation shown in Figure~\ref{fig:graph}. We consider $NB = 6$, generating a total of three iterations. At each iteration, we highlight the subset of batches to optimize. Nodes are labeled as: (1)~Optimized fixed batches (batches already optimized in previous iterations); (2)~Non-optimized fixed batches (batches not yet optimized but not selected in the given iteration); (3)~Batches to optimize (randomly selected batches to optimize in the given iteration); (4)~Unavailable batches (batches generated by the formulation but not available for the search procedure in the given iteration, according to variables ${MB}_k$).

Note that at the end of the search, each batch is optimized precisely once in one of the iterations. A higher diversification can be seen in this neighborhood, as it allows the selection of non-sequential batches. One can also note that in the last iteration (Iteration 3), only five batches are selected, although we defined $NB=6$. This is because these five batches are the only ones not optimized at this point in the search. The pseudo-code of the \textit{Multi-Batches Relocate} is shown in Algorithm~\ref{alg.multioprelocate}. To describe the procedure, we use $\cP$ to define the set of pairs machine/batch, where each element $(k, b) \in \cP$ represents a specific batch $b$ on a machine $k$, thus $\cP = \{(k, b) \,|\, k\in\cM, \, b \in~\cB_k\}$.

\renewcommand{\gets}{\leftarrow}
\begin{algorithm}[hbt!]
	\footnotesize
	\SetAlgoLined
	
	$NB \gets \lceil {\varphi} \times \sum_{k \in \cM}{MB_k} \rceil$;
	
	Initialize set $\cP$ of pairs ($k, b$) considering all machines $k \in \cM$ and their respective batches $b \in \cB_k \,|\, b \leq MB_k$;

	\While{$\cP \neq \varnothing$}{

		Select $NB$ pairs machine/batch $(k, b)$ randomly from $\cP$ to compose subset~$\cP'$;

		Create the subset~$\cO'$ of operations scheduled in the subset $\cP'$ of pairs machine/batch~$(k,b)$ on solution~$s$;

		Solve \textit{Batch Formulation}, for the given solution~$s$, for a subset of variables $X^{b}_{ik}$ in which $(k, b) \in \cP'$ and $i \in \cO'$;

		$\cP \gets \cP \setminus \cP'$;
		
		Update variables $MB_k$ according to Equation \ref{eq:mbcalc2}, if $\sum_{b \in \cB_k}{\eta^{b}_{k}}= MB_k$;
		
	}
	\caption{Multi-Batches Relocate ($s$, $\varphi$, $MB_k$)} 
	\label{alg.multioprelocate}
\end{algorithm}

Let $\cP' \subseteq \cP$ be the subset of selected pairs batch/machine on a given iteration. Algorithm~\ref{alg.multioprelocate} starts by computing the value of $NB$ (Line 1) and initializing the set $\cP$ (Line 2). The main loop of the algorithm (Lines 3-9) is repeated until there exist pairs batch/machine not optimized. Each iteration starts by randomly selecting $NB$ pairs batch/machine from the set $\cP$ to compose the subset $\cP'$ (Line 4). The subset $\cO'$ of operations assigned to any of the selected pairs batch/machine is created in Line 5. The \textit{Batch Formulation} is solved, for a limited time, for all variables $X^{b}_{ik}$, where $(k, b) \in \cP'$ and $i \in \cO'$ (Line 6). Then, the set $\cP$ and variables $MB_k$ are updated in Lines 7 and 8, respectively.

\section{Matheuristics}
\label{sec:matheuristics}

In this section, we present the matheuristics, which combines the \textit{MIP-based Neighborhood Searches} (Section~\ref{sec:mipneig}) and the \texttt{WMCT-WAVGA} heuristic (Algorithm~\ref{alg.constructive}) aiming to improve solutions continuously. Two well-known algorithm frameworks from the metaheuristics literature are considered, the Iterated Local Search (ILS) and the Greedy Randomized Adaptive Search Procedure (GRASP). The following sections explain each method that we refer to as \texttt{ILS-Math} and \texttt{GRASP-Math}, respectively.

For the local search, we consider a Variable Neighborhood Descent (VND) algorithm \citep{hansen2003variable}, using the two \textit{MIP-based Neighborhood Searches} described in Section~\ref{sec:mipneig}. First, we run the \textit{Multi-Batches Relocate} (Section \ref{sec:multibatchrelocate}), changing for the \textit{Batch Windows} (Section \ref{sec:batchwindows}) if no improvement is found. Every time an improved solution is found, we restart the local search, returning to the \textit{Multi-Batches Relocate}. The local search stops when no improvement is found after running both \textit{MIP-based Neighborhood Searches} completely. The sequence between the neighborhoods was defined based on preliminary experiments that showed a faster solution improvement using the \textit{Multi-Batches Relocate}.

\subsection{Iterated Local Search}
\label{sec:ils}

ILS is a powerful tool for combinatorial optimization problems with a simple structure and very useful for practical experiments. The method starts by building an initial solution and improving it using a local search. Then, the main loop consists of perturbing the current solution with simple modifications and running the local search until a stopping criterion is reached \citep{lourencco2003iterated}.

\renewcommand{\gets}{\leftarrow}
\begin{algorithm}[hbt!]
	\footnotesize
	\SetAlgoLined
	
	Build an initial PLSVSP solution $s_0$ using the \texttt{WMCT-WAVGA} heuristic (Algorithm~\ref{alg.constructive});
	
	Initialize variables $MB_k$ based on $s_0$ according to Equation \eqref{eq:mbcalc2};
	
	$s \gets VND (s_0, \rho, \varphi, MB_k)$;
	
	$s^{*} \gets s$;
	
	$\Omega \gets 1$;
	
	\While{$\Omega \leq \Omega^{max}$}{
		
		$\Omega \gets \Omega + 1$;
		
		$s' \gets RandomBatchSwap(s, \omega, MB_k)$;
		
		$s'^{*} \gets VND (s', \rho, \varphi, MB_k)$;
		
		\uIf{$f(s'^{*}) < f(s^{*}) \times (1 + \delta)$}{
			
			$s \gets s'^{*}$;
			
			\If{$f(s) < f(s^{*})$}{
				
				$s^{*} \gets s'^{*}$;
				
				$\Omega \gets 1$;
			}
		}
		\Else{
			$s \gets s^{*}$;
		}
		
	}
	
	Restart variables $MB_k$ based on the best solution $s^{*}$ according to Equation \eqref{eq:mbcalc2};
	
	$s^{*} \gets VND (s^{*}, \rho, \varphi, MB_k)$;
	
	\textbf{return} $s^*$;
	
	\caption{ILS-Math ($\rho$, $\varphi$, $\omega$, $\delta$,  $\Omega^{max}$)} 
	\label{alg.ils}
\end{algorithm}

The ILS matheuristic (\texttt{ILS-Math}) is described in Algorithm \ref{alg.ils}. A parameter $\Omega^{max} \in \mathbb Z_{+}$ defines the maximum number of iterations without improvement to execute, stopping the procedure whenever the value of the counter $\Omega$ reaches $\Omega^{max}$. During the ILS execution, worst solutions may be accepted, according to an acceptance parameter $\delta \in [0, 1]$. We use $s^{*}$ to keep the best solution found among all iterations.

The perturbation phase of the \texttt{ILS-Math} consists of randomly swapping batches in a given PLSVSP solution. In this approach, named \textit{RandomBatchSwap}, we compute the number of swaps ($NS$) to be performed based on a given parameter $\omega \in [0, 1]$, where $NS = \lceil {\omega} \times \sum_{k \in \cM}{MB_k} \rceil$. Note that the number of swaps to be performed is a fraction of the total number of available batches, considering all machines. At each swap movement, two batches are selected, and their operations are exchanged. If an empty batch is selected, the movement consists of removing the operations of the batch with operations inside and inserting them in the empty one. We forbid the selection of two empty batches. The movement is executed by updating the value of the corresponding assignment variables $X^{b}_{ik}$ of the selected batches.

The \texttt{ILS-Math} starts by building an initial solution $s_0$ (Line 1), using the \texttt{WMCT-WAVGA} heuristic (Algorithm \ref{alg.constructive}). In Line 2, the variable $MB_k$ is initialized (See Section \ref{sec:matheuristics}). The \textit{VND} is then performed on $s_0$ (Line 3), generating the current solution $s$, which is copied to $s^{*}$ (Line 4). The iteration counter $\Omega$ is initialized in Line 5. The main loop (Lines 6-19) is executed until $\Omega \leq \Omega^{max}$ and consists of repeatedly executing the \textit{Random Batch Swap} (Line 8), followed by the \textit{VND} (Line 9). In Line 10, the algorithm checks whether the objective value of the new solution $s'^{*}$, given by $f(s'^{*})$, passes the acceptance criteria. If true, the current solution $s$ is updated to $s'^{*}$ (Line 11), and the algorithm checks whether this solution is also better then $s^{*}$ (Line 12). If true, $s^{*}$ is updated (Line 13), and $\Omega$ is reinitialized (Line 14).  When the new solution is not accepted, the best solution $s^{*}$ replaces $s$  (Line 17). Finally, an intensification step is executed by running the \textit{VND} on the best solution $s^{*}$ (Lines 20-21).

\subsection{Greedy Randomized Adaptive Search Procedure}
\label{sec:GRASP}

GRASP is a multi-start method that combines a randomized constructive procedure followed by a local search, being successfully applied to many scheduling problems in the literature. For instance, we refer the reader to the papers of \cite{bassi2012planning}, \cite{rodriguez2012grasp}, and \cite{heath2013grasp}. In the constructive procedure, a Restricted Candidate List (RCL) with the most promising elements is built and one element is randomly selected at each step \citep{Resende2019}. To randomize the constructive procedure for the PLSVSP, we replace the operation selection step of the \texttt{WMCT-WAVGA} heuristic (Line 6 of Algorithm \ref{alg.constructive}) by a \textit{Randomized Operation Selection}, using a parameter $\alpha \in [0, 1]$ to define the greediness of the method. When $\alpha = 0$, the method builds the same solution of the \texttt{WMCT-WAVGA} heuristic, and when $\alpha = 1$, a completely randomized solution is generated. After computing the priority value $\pi_i$ for the set $\cU$ of unscheduled operations (Line 5 of Algorithm \ref{alg.constructive}), we identify the minimum and maximum priority values, defined as $\pi_{min}$ and $\pi_{max}$, respectively. The RCL is created as RCL $= \{i \in \mathcal{U} \:|\: \pi_i \geq \pi_{max} - \alpha (\pi_{max}-\pi_{min})\}$, and one operation $i^*$ in randomly selected from the RCL.

\renewcommand{\gets}{\leftarrow}
\begin{algorithm}[hbt!]
	\footnotesize
	\SetAlgoLined
	
	$\Omega \gets 1$;
	
	$s^{*} \gets \emptyset$; $f(s^{*}) \gets \infty$;
	
	\While{$\Omega \leq \Omega^{max}$}{
		
		$\Omega \gets \Omega + 1$;
		
		Build a PLSVSP solution $s$ using the \textit{Randomized Constructive Procedure};
		
		Initialize variables $MB_k$ based on $s$ according to Equation \eqref{eq:mbcalc2};
		
		$s'^{*} \gets VND (s, \rho, \varphi, MB_k)$;
		
		\If{$f(s'^{*}) < f(s^{*})$}{
			$s^{*} \gets s'^{*}$;
			
			$\Omega \gets 1$;
		}
		
	}
	
	Restart variables $MB_k$ based on the best solution $s^{*}$ according to Equation \eqref{eq:mbcalc2};
	
	$s^{*} \gets VND (s^{*}, \rho, \varphi, MB_k)$;
	
	\textbf{return} $s^*$;
	
	\caption{GRASP-Math ($\rho$, $\varphi$, $\alpha$,  $\Omega^{max}$)} 
	\label{alg.grasp}
\end{algorithm}

The pseudo-code of the GRASP matheuristic (\texttt{GRASP-Math}) is shown in Algorithm \ref{alg.grasp}. We call \textit{Randomized Constructive Procedure}, the \texttt{WMCT-WAVGA} heuristic (Algorithm \ref{alg.constructive}) with the \textit{Randomized Operation Selection}. A parameter $\Omega^{max} \in \mathbb Z_{+}$ defines the maximum number of iterations without improvement, while $s^{*}$ saves the best solution found among all iterations.

The \texttt{GRASP-Math} starts by initializing the $\Omega$ counter and $f(s^*)$ (Lines 1 and 2). The main loop is executed until $\Omega \leq \Omega^{max}$ (Lines 3-12). At each iteration, a new randomized solution $s$ is built (Line 5), with its $MB_k$ variables initialized (Line 6), and the \textit{VND} is applied to this solution (Line 7), generating a new solution $s'^{*}$. The algorithm checks whether the objective value of the new solution, given by $f(s'^{*})$, is better than $f(s^*)$ (Line 8). If true, $s^{*}$ is replaced by $s'^{*}$, and the counter $\Omega$ is reinitialized (Lines 9-10). Finally, an intensification step is executed by running the \textit{VND} on the best solution $s^*$ (Lines 13-14).

\subsection{Overview of the Methodology}
\label{sec:general_math}

The general pseudocode of the proposed methodology is described in Algorithm \ref{alg.general}. The algorithm allows selecting different mathematical formulations to be used in the matheuristics' main loop and the intensification step.

\renewcommand{\gets}{\leftarrow}
\begin{algorithm}[hbt]
	\footnotesize
	\SetAlgoLined
	
	Choose a matheuristic framework between \texttt{GRASP-Math} and \texttt{ILS-Math}.
	
	Select a mathematical formulation between \texttt{Batch-WSPT} and \texttt{Batch-S}.
	
	Run the chosen matheuristic's main loop  using the selected mathematical formulation. 
	
	Choose a mathematical formulation between \texttt{Batch-WSPT} and \texttt{Batch-S}.
	
	Initialize the variables of the chosen mathematical formulation according to the best solution found so far.
	
	Run the VND algorithm using the chosen mathematical formulation.
	
	Return the best-found solution.
	
	\caption{General Methodology} 
	\label{alg.general}
\end{algorithm}

The combination of the developed matheuristics' frameworks (\texttt{ILS\--Math} and \texttt{GRASP\--Math}) and the batch scheduling formulations (\texttt{Batch\--WSPT} and \texttt{Batch-S}), following the defined methodology, generates three variants of each matheuristic, which we refer to as \MHa, \MHb, \MHc, \MHd, \MHe, and \MHf, described below:

\begin{itemize}
	\item \MHa: \: \:\:\: \texttt{ILS-Math} with \texttt{Batch-WSPT}
	\item \MHb: \: \:\:\: \texttt{ILS-Math} with \texttt{Batch-S}
	\item \MHc: \: \:\:\: \texttt{ILS-Math} with \texttt{Batch-WSPT} and \texttt{Batch-S}
	\item \MHd: \: \texttt{GRASP-Math} with \texttt{Batch-WSPT}
	\item \MHe: \: \texttt{GRASP-Math} with \texttt{Batch-S}
	\item \MHf: \: \texttt{GRASP-Math} with \texttt{Batch-WSPT} and \texttt{Batch-S}
\end{itemize}

The \MHc\: and the \MHf\: consider the \texttt{Batch-WSPT} formulation for their main loop, changing for the \texttt{Batch-S} formulation at the intensification step.

\section{Computational Experiments}
\label{sec:experiments}

In this section, we present the computational experiments conducted to assess the performance of the proposed matheuristics. We compare them with the \textit{Mathematical Formulations} (Section \ref{sec:formulation}) running independently. The computational experiments were performed on a machine with an Intel i7-8700K CPU of 3.70GHz and 64 GB of RAM running Linux. All methods were coded using \Cpp language solved by CPLEX 12.8 solver running in a single thread with MIP emphasis set to finding hidden feasible solutions. We limit the execution time to one second per CPLEX call in the matheuristics. Thus, sub-problem optimizations are interrupted within the defined time-limit even if the optimal solutions have not been reached.

To evaluate the solution's quality, we compare the solutions provided by each method with the Best-Know Solutions (BKS) from the literature in terms of the Relative Percentage Deviation (RPD), computed according to Equation \eqref{eq:RPD}. $TWCT^{Method}$ designates the total weighted completion time obtained with one run of a selected method on a PLSVSP instance, while $TWCT^{BKS}$ is the total weighted completion time for the Best-Know Solution for the same instance. In all analyzes, we compare the matheuristics, running each one ten times per instance, with the solutions of \texttt{Batch-WSPT} (Section \ref{sec:formulation1}) and \texttt{Batch-S} (Section \ref{sec:formulation2}) formulations, running within a 6-hour time-limit. We limit the memory allocation to 10 GB for each method execution to allow multiple runs simultaneously using the available CPUs.

{\small
\begin{flalign}
    \label{eq:RPD}
    RPD = 100 \times \frac{\mathit{TWCT}^{Method} - \mathit{TWCT}^{\mathit{BKS}}}{ \mathit{TWCT}^{\mathit{BKS}}}
\end{flalign}
}

The discussion of the results focuses only on comparing the \MHc\: and the \MHf\: matheuristics with the pure mathematical formulations since these approaches presented the best results among the proposed matheuristic variants in a preliminary analysis, shown in Appendix \ref{app:mathcomp}. We also included the complete results for each instance, considering all methods, in Appendix \ref{app:comp_results}.

\subsection{Instances Description}
\label{sec:inst_desc}

The experiments were conducted on the benchmark of 72 PLSVSP instances\footnote{The PLSVSP instances benchmark is available online \cite{abumarrul2019}}, developed by \cite{AbuMArrul2020} in which the number of operations to schedule varies as $\abs{\cO} \in \{15, 25, 50\}$, and the number of machines varies as $|\cM| \in \{4, 8\}$. The number of jobs is defined as $|\cN|=\lfloor|\cO|/3\rfloor$, and the number of families is defined as $|\cF|=3$. The instances generation process, defined by the authors, with the distributions used, is detailed below.

Operations are randomly assigned to families, and family setup times are drawn from a discrete uniform distribution $U(5,10)$. Processing times of operations follow a discrete uniform distribution $U(1,30)$, and their occupancy values are drawn from another discrete uniform distribution $U(1,100)$ with step size 10. Job's weights follows a discrete uniform distribution $U(1,50)$. Capacities of machines are drawn from a discrete uniform distribution $U(80,100)$ with step size 10. Release dates of machines and operations follow a discrete uniform distribution $U(0,MR)$, where $MR$ is the maximum release date, defined as $MR = \bigl\lceil \alpha \cdot  \sum_{i \in \cO}{\bigl(p_i + s_{f_i}\bigr)}/ |\cM| \bigr\rceil$. The release factor $\alpha \in \{0.25, 0.5, 0.75\}$ defines different values for $MR$. An eligibility factor $\beta \in \{0.7, 0.9\}$, defines the probability of a machine being eligible to process each operation. Finally, a factor $\gamma \in \{0.05, 0.15\}$ defines the probability of a job being associated with an operation.

\subsection{Parameter Tuning}
\label{sec:partun}

For parameter tuning, we used 12 medium-sized instances (25 operations) and 12 large-sized instances (50 operations), selected at random, corresponding to one-third of the total number of instances available in the benchmark. The \texttt{Batch-WSPT} formulation was used within the matheuristics during the parameterization. First, we set parameters related to the \textit{MIP-based Neighborhood Searches} (Section \ref{sec:mipneig}), to further define each specific matheuristic parameter. To establish the neighborhoods parameters, we ran five times each neighborhood individually, using the solution generated by the \texttt{WMCT-WAVGA} heuristic (Algorithm \ref{alg.constructive}) as a warm start. We define the ranges for the \textit{Batch Windows} (Section \ref{sec:batchwindows}), and for the \textit{Multi-Batches Relocate} (Section \ref{sec:multibatchrelocate}) parameters as $\rho \in [0.1, 0.5]$, and $\varphi \in [0.1, 0.5]$, respectively, considering a step of size $0.05$. After setting values for $\rho$, and $\varphi$, the \texttt{ILS-Math} was executed five times without accepting worse solutions. We used $\Omega^{max} = 10$, limiting $\omega \in [0.05, 0.15]$ with a step size of $0.05$, to define the perturbation parameter. With $\omega$ set, we ran the \texttt{ILS-Math} again five more times, limiting $\delta \in [0.00, 0.15]$ with a step size of $0.05$, to define the worst solution acceptance rate. The same steps were executed to define the greediness factor of the \texttt{GRASP-Math}, limiting $\alpha \in [0.05, 0.15]$, with a step size of $0.05$. The final values for the parameters are shown in Table \ref{tab:parametersdef}.

\begin{center}
	\begin{table}[htbp]
		\centering
		\footnotesize
		\caption{Parameters definition}	
		\begin{tabularx}{\textwidth}{lcXcc}
		\hline
			Algorithm & Parameter & Description & Domain & Value \\
             \hline
             Batch Windows & $\rho$ & The proportion of the makespan ($C_{max})$ to optimize at each iteration & [0, 1] & 0.20 \\
               Multi-Batches Relocate & $\varphi$ & The proportion of batches to optimize at each iteration & [0, 1] & 0.30 \\
                Random Batch Swap & $\omega$ & The proportion of batch swaps to execute at each perturbation & [0, 1] & 0.10 \\
                ILS-Math & $\delta$ & Worse solutions acceptance rate & [0, 1] & 0.00 \\    
                GRASP-Math & $\alpha$ & Greediness factor for the randomized constructive procedure & [0, 1] & 0.10 \\
                 ILS-Math and GRASP-Math & $\Omega^{max}$ & Maximum number of iterations without improvement & $\mathbb{Z}_{+}$ & 10 \\
			\bottomrule			
		\end{tabularx}%
		
		\label{tab:parametersdef}%
	\end{table}%
\end{center}

\subsection{Results Analysis and Discussion}
\label{sec:ilsmath}

\subsubsection{Average Values for the RPD and Computational Time}
\label{sec:avgvalue}

In the first analysis, shown in Table \ref{tab:mean_results}, we compare the four methods (\texttt{Batch-WSPT}, \texttt{Batch-S}, \MHc\, and \MHf\,) in terms of the average RPD ($\overline{RPD}$) and the average computational time ($\overline{time}$), in seconds, for each instance group. We use the grouping scheme defined by \cite{AbuMArrul2020} for this analysis, in which each group, represented by the combination of the number of operations and machines ($|\cO|$ and $|\cM|$), comprises 12 instances. The best result for each criterion in each group is highlighted in bold.

\begin{table}[htbp]
 \centering
 \footnotesize
  \caption{Average values for the RPD and computational time distributions for each method in each instance group.}
  \resizebox{\textwidth}{!}{%
  \begin{threeparttable}
    \begin{tabular}{ccrccccccccccc}
    \toprule
    \multirow{2}[4]{*}{$|\cO|$} & \multirow{2}[4]{*}{$|\cM|$} &       & \multicolumn{3}{c}{\texttt{Batch-WSPT}} & \multicolumn{2}{c}{\texttt{Batch-S}} &       & \multicolumn{2}{c}{\MHc} &       & \multicolumn{2}{c}{\MHf} \\
\cmidrule{4-8}\cmidrule{10-11}\cmidrule{13-14}          &       &       & $\overline{RPD}$   & $\overline{time}$   &       & $\overline{RPD}$   & $\overline{time}$   &       & $\overline{RPD}$   & $\overline{time}$   &       & $\overline{RPD}$   & $\overline{time}$ \\
    \midrule
    \multirow{2}[1]{*}{15} & 4     &       & 0.49  & 20414 &       & \textbf{0.00} & 20725 &       & 0.42  & 14    &       & 0.30  & 26 \\
          & 8     &       & 0.04  & 18000 &       & \textbf{0.01} & 18000 &       & 0.08  & 46    &       & \textbf{0.01} & 54 \\
          &       &       &       &       &       &       &       &       &       &       &       &       &  \\
    \multirow{2}[0]{*}{25} & 4     &       & 0.70  & 21600 &       & 0.81  & 21600 &       & -0.17 & \textbf{123} &       & \textbf{-0.31} & 186 \\
          & 8     &       & 0.23  & 21600 &       & 0.62  & 21600 &       & -0.20 & \textbf{244} &       & \textbf{-0.27} & 297 \\
          &       &       &       &       &       &       &       &       &       &       &       &       &  \\
    \multirow{2}[1]{*}{50} & 4     &       & -0.27 & 19978$^{\dagger}$ &       & 0.60  & 21600 &       & \textbf{-3.02} & \textbf{1052} &       & -2.65 & 1159 \\
          & 8     &       & -1.83 & 17380$^{\dagger}$ &       & 1.90  & 13901$^{\dagger}$ &       & \textbf{-4.09} & \textbf{990} &       & -3.80 & 1066 \\
    \midrule
    \multicolumn{2}{c}{\multirow{2}[2]{*}{All instances}} &       &       &       &       &       &       &       &       &       &       &       &  \\
    \multicolumn{2}{c}{} &       & -0.11 & 19829$^{\dagger}$ &       & 0.66  & 19571$^{\dagger}$ &       & \textbf{-1.16} & \textbf{412} &       & -1.12 & 465 \\
    \bottomrule
    \end{tabular}%
    \begin{tablenotes}
        \item[$\dagger$] CPLEX execution interrupted before 21,600 seconds (time limit) for some instances in this group, for this method, due to memory problems.
    \end{tablenotes}
    \end{threeparttable}
    
    }
  \label{tab:mean_results}%
\end{table}%

Note that the average relative percentage deviations are low for instances with 15 operations, with an $\overline{RPD}$ below 1\% for all methods. Regarding group 15--4, the \texttt{Batch-S} formulation presented the smallest $\overline{RPD}$ with $0.00$, that is, it reaches the BKS in all runs, but with an average computational time of $20{,}725$ seconds. The average computational times for the MIP formulations in this group are close to the limit of $21{,}600$ seconds (6 hours) previously defined. \MHc\, runs with the least average computational time in this group (14 seconds). However, concerning the $\overline{RPD}$, \MHf\, is the best matheuristic approach ($\overline{RPD}=0.30$), without significantly increasing computational time spent by \MHc. Similar behavior can be seen in group 15--8. \MHf\, presented equivalent results with \texttt{Batch-S} ($\overline{RPD} = 0.01$), but consuming less than a minute to achieve it. In contrast, \texttt{Batch-S} run for an average of $18,000$ seconds.

For medium-sized groups (25--4 and 25--8), the matheuristics outperform the MIP formulations in terms of solution quality, with considerably less computing time. The MIP formulations maintained a good quality of the solutions, below 1\% for the $\overline{RPD}$ in both groups, but running until the time limit in all executions ($\overline{time} = 21{,}600$). Between the matheuristics, \MHf\, outperforms \MHc\, in terms of solution quality, with $\overline {RPD} = -0.31$ in group 25--4, and $\overline {RPD} = -0.27$ in group 25--8, but with 51\% extra computational time needed on average in the former group ($186$ seconds against $123$ seconds), and 22\% in the latter ($297$ seconds against $244$ seconds).

For large-sized groups, we can see that \texttt{Batch-S} loses performance, with the worst values for the $\overline{RPD}$ in both groups ($0.60$ in group 50--4 and $1.90$ in group 50--8), resulting in a clear advantage for \texttt{Batch-WSPT}. Again, the matheuristics dominates the MIP formulations in both groups with lower values for the $\overline{RPD}$ and $\overline{time}$. However, in these groups, \MHc\, dominates \MHf\,, in terms of solution quality and computational time. It can be noted that the $\overline{time}$ for the MIP formulations on these groups is smaller than the limit of $21{,}600$. This is due to memory issues, with the solver interrupting the execution of 6 instances, in each formulation, before reaching the time limit. Appendix \ref{app:comp_results} details the complete results, indicating these instances.

One can note that, in all groups, the GRASP approach requires more computational time than ILS due to the method's restart characteristic. This characteristic also affects the quality of the method's solutions when the size of the instances increases. We can also observe that the average computational time for the matheuristics grows considerably when the number of operations to schedule increases. It is known that mathematical models are severely affected when the size of problems grows, which, in our case, directly impacts the execution time of our matheuristics. 

\subsubsection{RPD Distribution for Different Instances Aspects}
\label{sec:rpd_dist}

\begin{figure}[htb!]
	\subfigure[15 operations.\label{fig:boxplot15op}]{%
		\includegraphics[width=0.315\textwidth, trim = 10 28 10 0, clip]{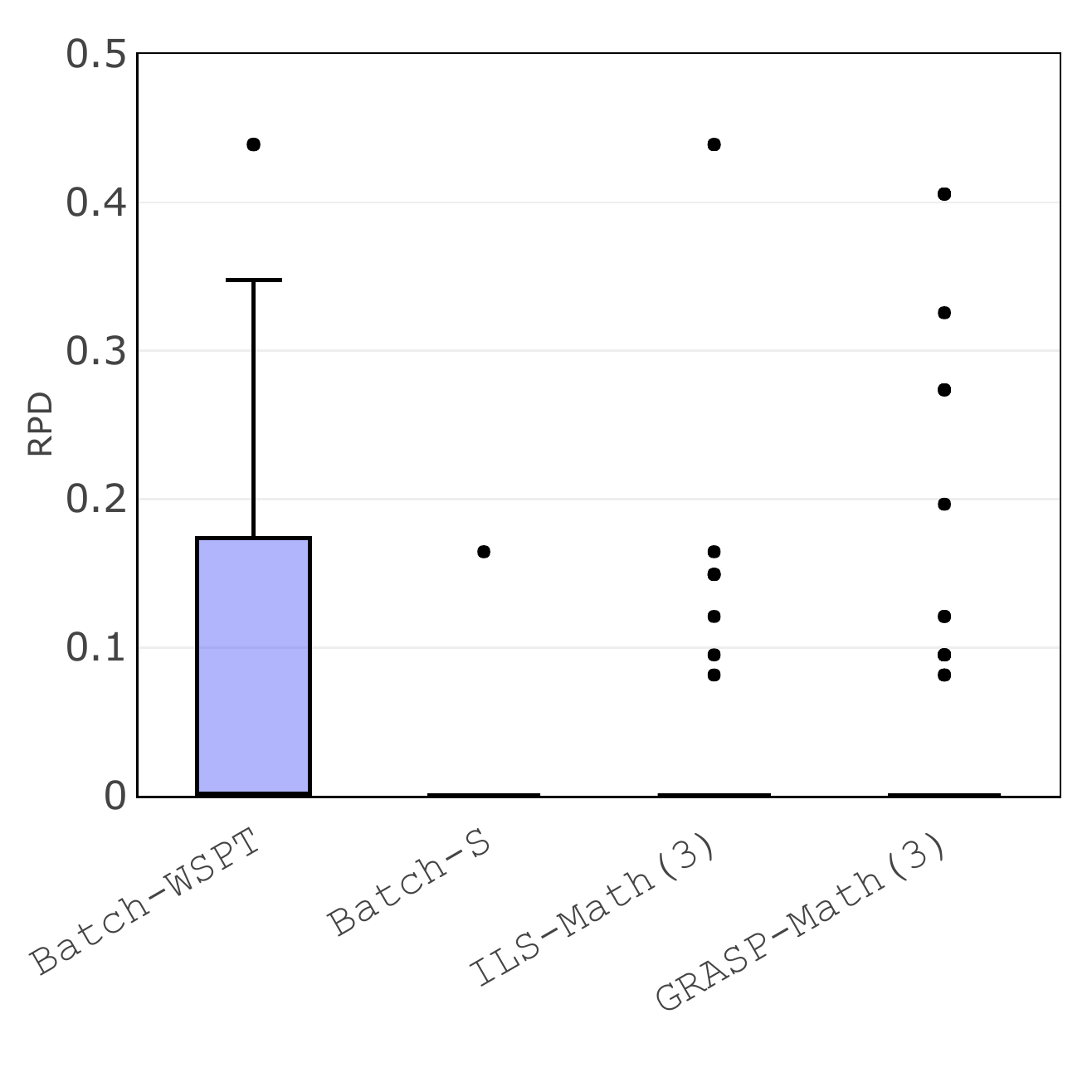}
	}
	\subfigure[25 operations.\label{fig:boxplot25op}]{%
		\includegraphics[width=0.315\textwidth, trim = 10 28 10 0, clip]{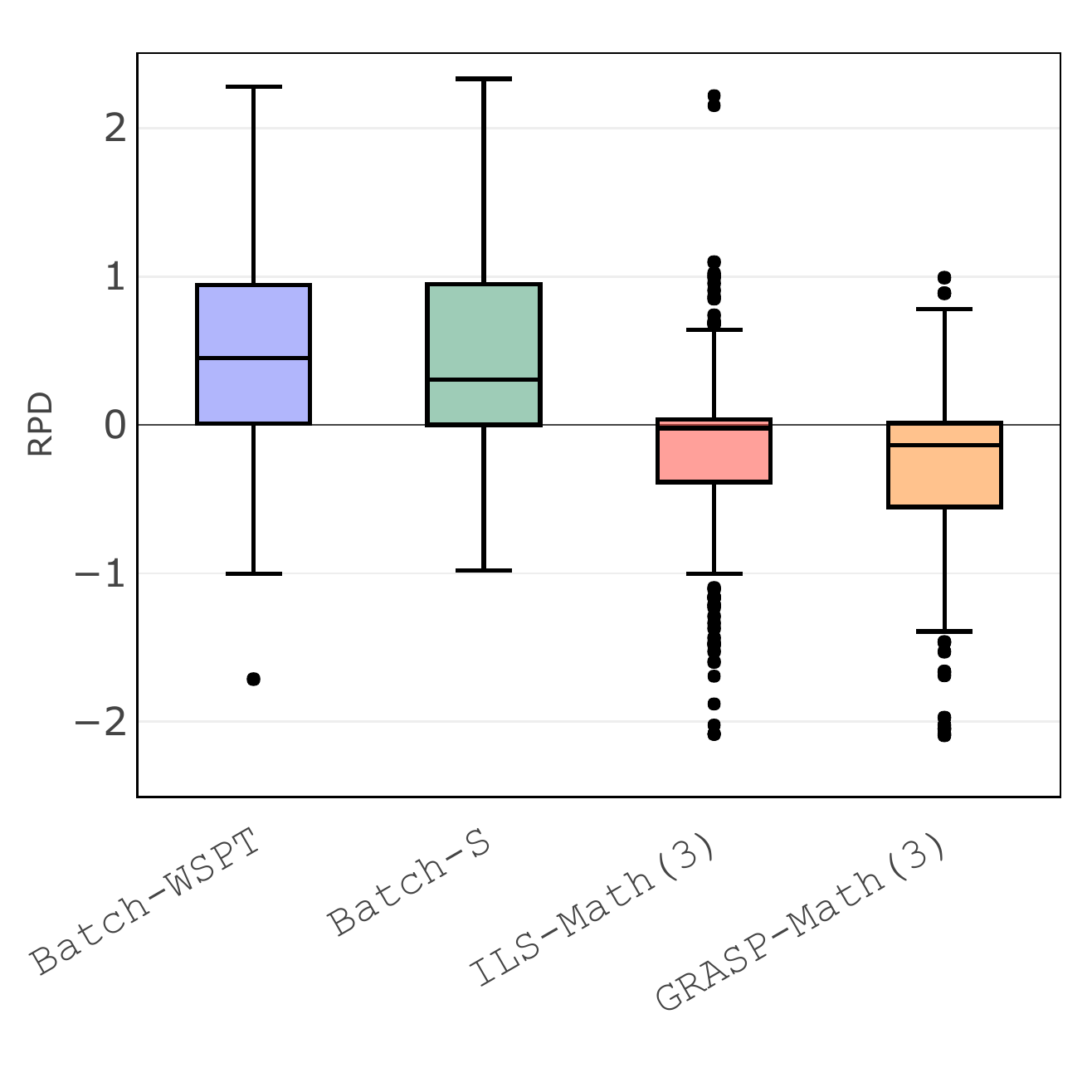}
	}
	\subfigure[50 operations.\label{fig:boxplot50op}]{%
		\includegraphics[width=0.315\textwidth, trim = 10 28 10 0, clip]{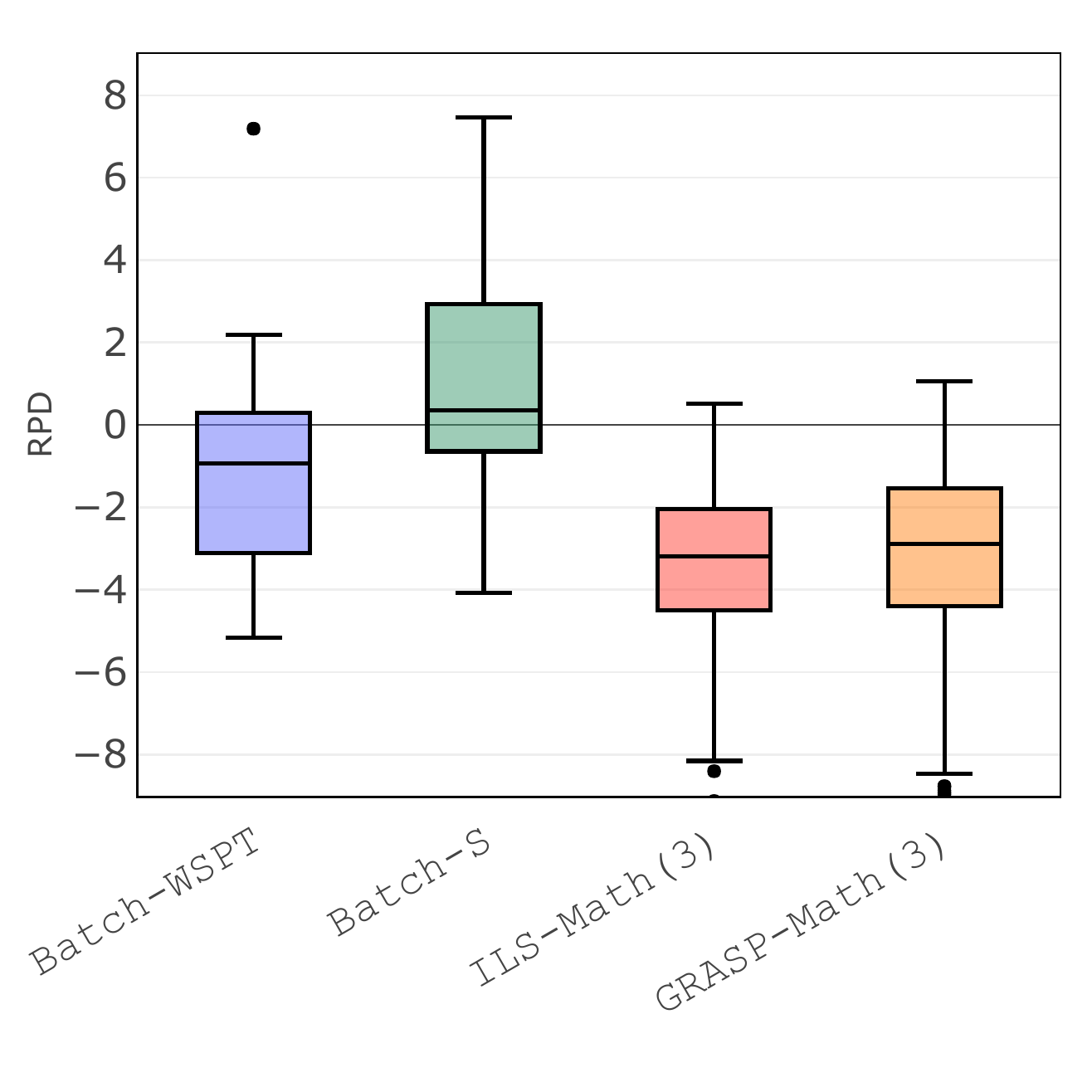}
	}
	\caption{Boxplots of the RPD distributions for each method, considering different number of operations.}
	\label{fig:boxplot_operations}
\end{figure}

\begin{figure}[http!]
\centering
	\subfigure[4 machines.\label{fig:boxplot4mach}]{%
		\includegraphics[width=0.325\textwidth, trim = 10 28 10 0, clip]{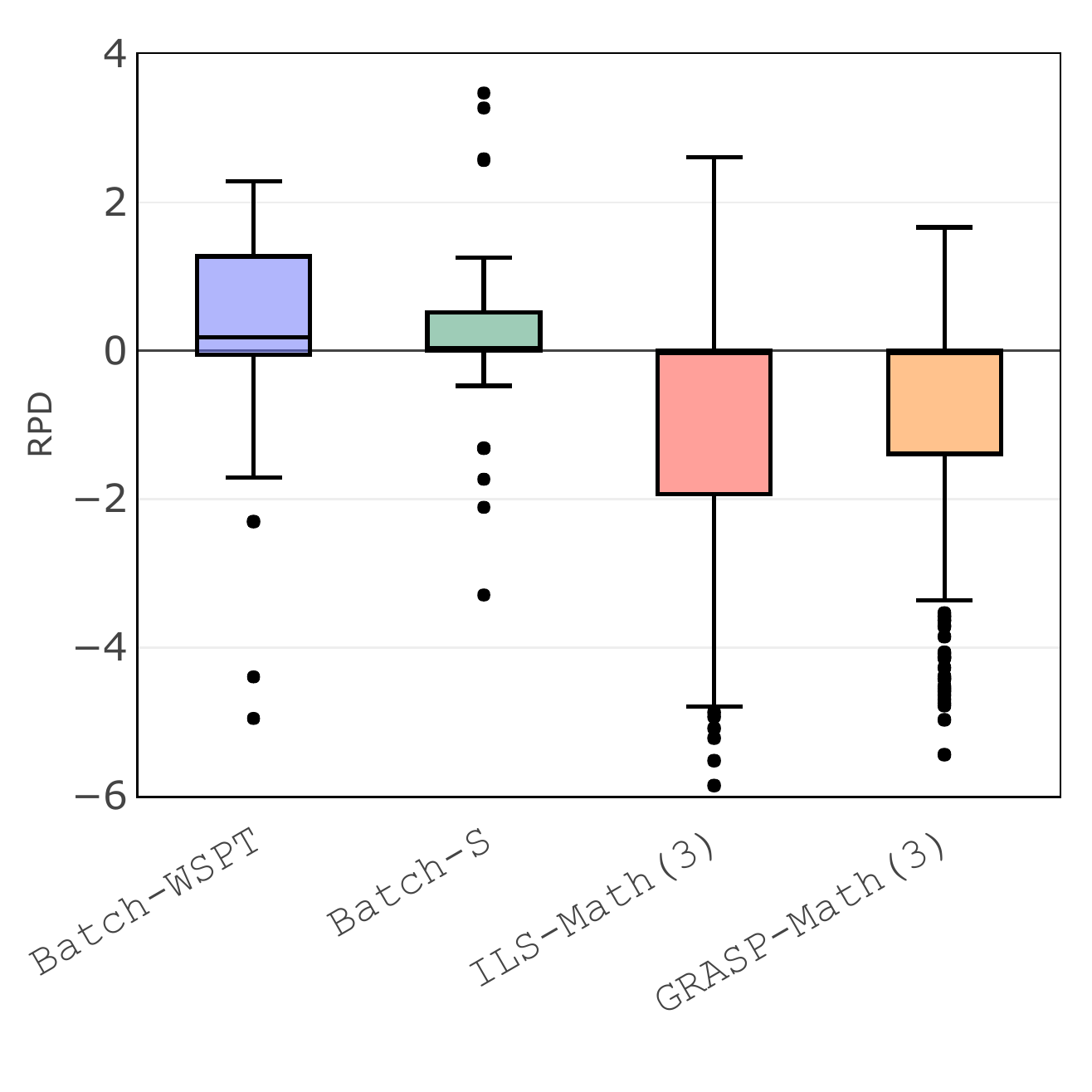}
	}
	\subfigure[8 machines.\label{fig:boxplot8mach}]{%
		\includegraphics[width=0.325\textwidth, trim = 10 28 10 0, clip]{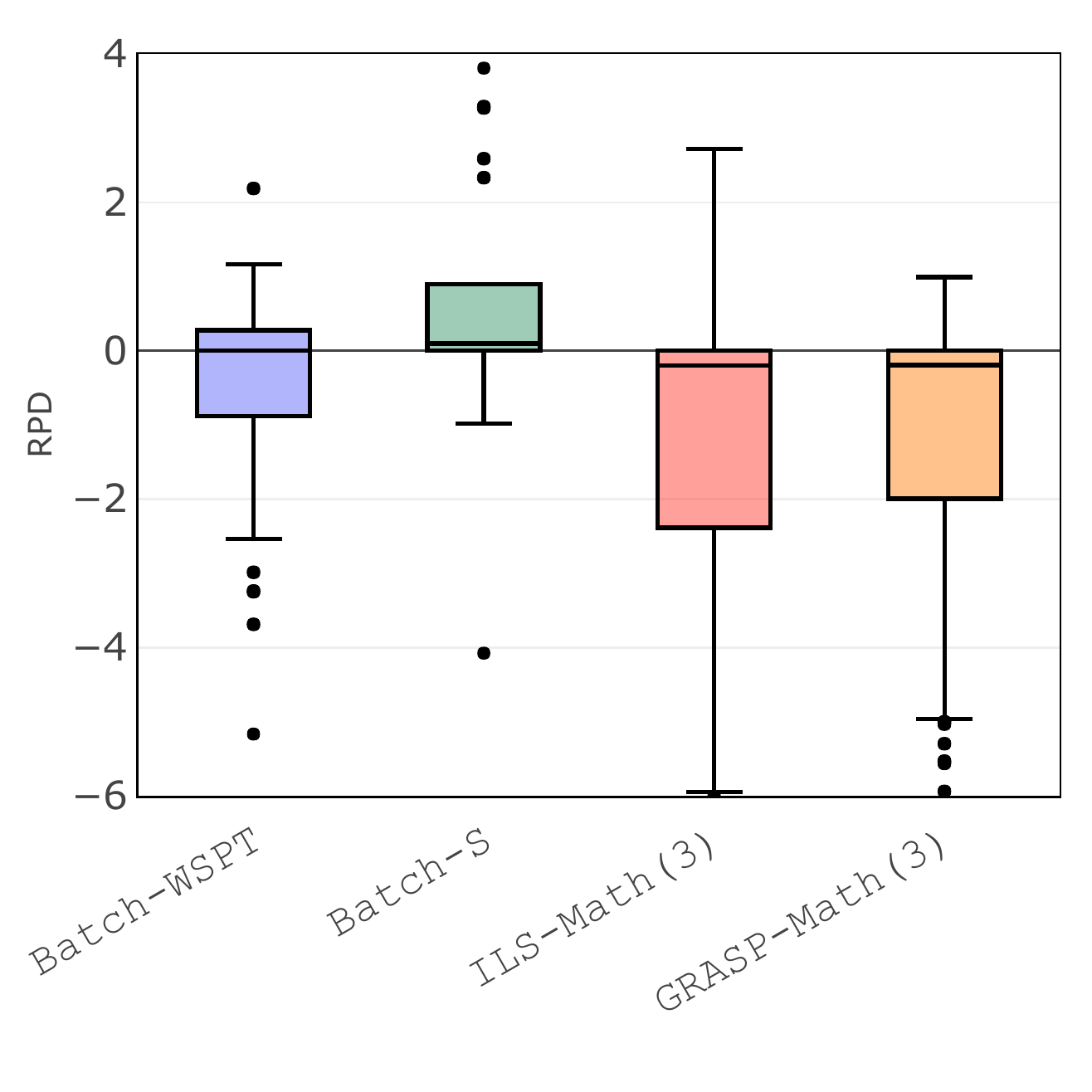}
	}
	\caption{Boxplots of the RPD distributions for each method, considering different number of machines.}
	\label{fig:boxplot_machines}
\end{figure}

In this analysis, we evaluate the impacts on the solution quality of each method with the increasing number of machines and operations. Figures \ref{fig:boxplot15op}, \ref{fig:boxplot25op}, and \ref{fig:boxplot50op} depict the distribution of the RPDs for instances with 15, 25, and 50 operations to schedule, respectively. Next, Figures \ref{fig:boxplot4mach}, and \ref{fig:boxplot8mach}, show the distribution of the RPDs for instances with 4, and 8 machines, respectively.

Regarding instances with 15 operations (Figure \ref{fig:boxplot15op}), all methods are competitive (RPDs closest to zero). However, \texttt{Batch-WSPT} struggles to achieve the BKS due to the reduced solution space considered by the formulation. Note that the matheuristics have better performance than the MIP formulations when 25 and 50 operations are considered (Figures \ref{fig:boxplot25op} and \ref{fig:boxplot50op}). It is interesting to see the distributions' behavior when the number of operations increases. For the MIP formulations, we can see that \texttt{Batch-S} performs better on small instances, while \texttt{Batch-WSPT} shows advantages when the number of operations increases. Note that on medium-sized instances (25 operations), they have similar performance. The consideration of extra variables and constraints to sequence operations within batches proved to be a good strategy for small instances since it considers the complete solution space of the problem, however affecting the performance of \texttt{Batch-S} when the number of operations increases. Regarding the matheuristics, we can see that they maintain similar distributions regardless of the number of operations, but with \MHc\, improving its performance, compared to \MHf\, in instances with 50 operations.

The same behavior can be observed in the analysis by the number of machines (Figures \ref{fig:boxplot4mach} and \ref{fig:boxplot8mach}), with the matheuristics showing better distribution of RPDs. However, we can observe less dispersed distributions due to the consideration of instances with 15 operations in both groups. The MIP formulations maintained the same behavior, with Batch-S losing performance when the number of machines grows.

\subsubsection{Comparison with the BKS Values}
\label{sec:bks_comp}

Now, we analyze the methods, in Table \ref{tab:bks_analysis}, regarding the achievement or improvement of the Best-Know Solutions from the literature, defined by \cite{AbuMArrul2020}. The results indicate the number of instances in which the BKS was achieved or improved (\#Inst.), its percentage regarding the complete set of 72 PLSVSP instances (\%Inst.), and the respective percentage of runs of achievement or improvement (\%Runs). The $\overline{RPD}$ for the subsets of not achieved and not improved solutions are included. Finally, we show the minimum RPD ($RPD^-$) found by each method when the BKS is improved.

\begin{table}[htbp]
  \centering
  \footnotesize
   \setlength{\tabcolsep}{1.5pt}
  \caption{Analysis regarding BKS solutions from the literature.}
  \resizebox{\textwidth}{!}{%
  \begin{threeparttable}
    \begin{tabular}{lrccccccccccccc}
    \toprule
    \multirow{2}[4]{*}{Method} &       & \multicolumn{3}{l}{Achieved BKS$^{\dagger}$} &       & Not Achieved &       & \multicolumn{3}{l}{Improved BKS} &       &       & \multicolumn{1}{l}{Not Improved} \\
\cmidrule{3-5}\cmidrule{7-7}\cmidrule{9-12}\cmidrule{14-14}          &       & \#Inst. & \%Inst. & \%Runs &       & $\overline{RPD}$   &       & \#Inst. & \%Inst. & \%Runs & $RPD^-$ &       & $\overline{RPD}$ \\
    \midrule
    \texttt{Batch-WSPT} &       & 42    & 58.33 & 58.33 &       & 1.16  &       & 23    & 31.94 & 31.94 & -5.17 &       & 0.71 \\
    \texttt{Batch-S} &       & 41    & 56.94 & 56.94 &       & 2.06  &       & 13    & 18.06 & 18.06 & -4.08 &       & 1.08 \\
    \MHc &       & 69    & 95.83 & 81.94 &       & 0.72  &       & 40    & 55.56 & 49.72 & -10.96 &       & 0.26 \\
    \MHf &       & 72    & 100.00 & 86.11 &       & 0.59  &       & 40    & 55.56 & 50.56 & -9.72 &       & 0.17 \\
    \bottomrule
    \end{tabular}%
        \begin{tablenotes}
        \item[$\dagger$] This subset includes solutions that have also improved the BKS. 
    \end{tablenotes}
    \end{threeparttable}
    
    }
  \label{tab:bks_analysis}%
\end{table}%

Note that all methods achieve the BKS in more than 50\% of instances, with \MHf\, reaching it in all instances (\#Inst. = 72). The \%Inst. and \%Run has the same value for each MIP formulation, as we only run it once for each instance. Again, the matheuristics have better performance than the MIP formulations. Note that the $\overline{RPD}$ for runs where BKS was not found is less than 1\% for the matheuristics and more than 1\% for the MIP formulations, with the worst performance for \texttt{Batch-S} formulation ($\overline{RPD}=2.06\%$). It can be seen that the matheuristics \MHc\ and \MHf\, achieve the BKS in 81.94\% and 86.11\% of the executions, respectively. The same advantage for the matheuristics can be seen for improved solutions, with both methods improving the BKS in about 50\% of runs. The $\overline{RPD}$ is low for runs that have not improved the BKS (0.26\% for \MHc\,, and 0.17\% for \MHf), as it includes runs that have achieved the BKS. With these results, we can see that 40 new best solutions are found by the matheuristics, with an improvement of 10.96\% in the objective value of the solution in the best case ($RPD^-$ of \MHc\,). The specific instances with enhanced solutions are indicated in Appendix \ref{app:comp_results}. Considering that 26 instances have the proven optimal solutions (see Table \ref{tab:comp_results}), 40 new best solutions have been defined in 46 possible instances.

\subsubsection{Statistical Analysis for the RPD Distributions}
\label{sec:wcx_comp}

To improve our discussion and validate what we highlighted during the previous analyses, we applied a statistical evaluation comparing the RPD distributions between the methods, considering the complete set of 72 PLSVSP instances. First, we tested the distributions' normality with the Shapiro-Wilk test, which shows that the RPDs do not follow a normal distribution. Then, we ran the pairwise Wilcoxon rank-sum test with Hommel's $p$-values adjustment. We also included the Analysis of Variance (ANOVA) with the Tukey HSD (honestly significant difference) test to compare the method's RPD distributions. We ran both tests with a confidence level of 0.05. In Table \ref{tab:wilcoxom}, we present the $p$-value for each pair of methods, also including some statistics on the RPD distributions. Significantly better results are highlighted in bold.

\begin{table}[htbp]
	\centering
	\footnotesize
	\caption{The mean and standard deviation of the RPD distributions for each method, and $p$-values from pairwise Wilcoxon rank-sum and ANOVA with Tukey HSD tests, with a 0.05 confidence level, considering all 72 PLSVSP instances.}
	\begin{tabular}{lcccccc}
		\toprule
		\multicolumn{1}{c}{\multirow{2}[4]{*}{Method}} & \multirow{2}[4]{*}{$\overline{RPD}$} & \multirow{2}[4]{*}{$\sigma$} &       & \multicolumn{3}{c}{Wilcoxon $p$-value (ANOVA--Tukey $p$-value)} \\
		\cmidrule{5-7}          &       &       &       & \texttt{Batch-WSPT} & \texttt{Batch-S} & \MHf \\
		\midrule
		\texttt{Batch-WSPT} & -0.10 & 1.78  &       &       &       &  \\
		\texttt{Batch-S} & 0.66  & 2.10   &       & 0.38 (0.14)  &       &  \\
		\MHf & -1.16 & 2.18  &       & \textbf{0.00 (0.00)} & \textbf{0.00 (0.00)} &  \\
		\MHc & -1.12 & 2.15  &       & \textbf{0.00 (0.00)} & \textbf{0.00 (0.00)} & 0.50 (0.98) \\
		\bottomrule
	\end{tabular}%
	
	\label{tab:wilcoxom}%
\end{table}%

The tests confirm that the matheuristics yield significantly better results than the MIP formulations. Note that no statistical significance can be seen when comparing one MIP formulation with another, nor when comparing the matheuristics between each other. Wilcoxon and ANOVA--Tukey agreed in all cases despite discrepancies in the $p$-values.

\subsubsection{Average RPD Evolution Analysis}
\label{sec:rpd_evol}

\begin{figure} [htbp!] 
	\centerline{\includegraphics[scale=0.48, trim = 0 0 0 0, clip]{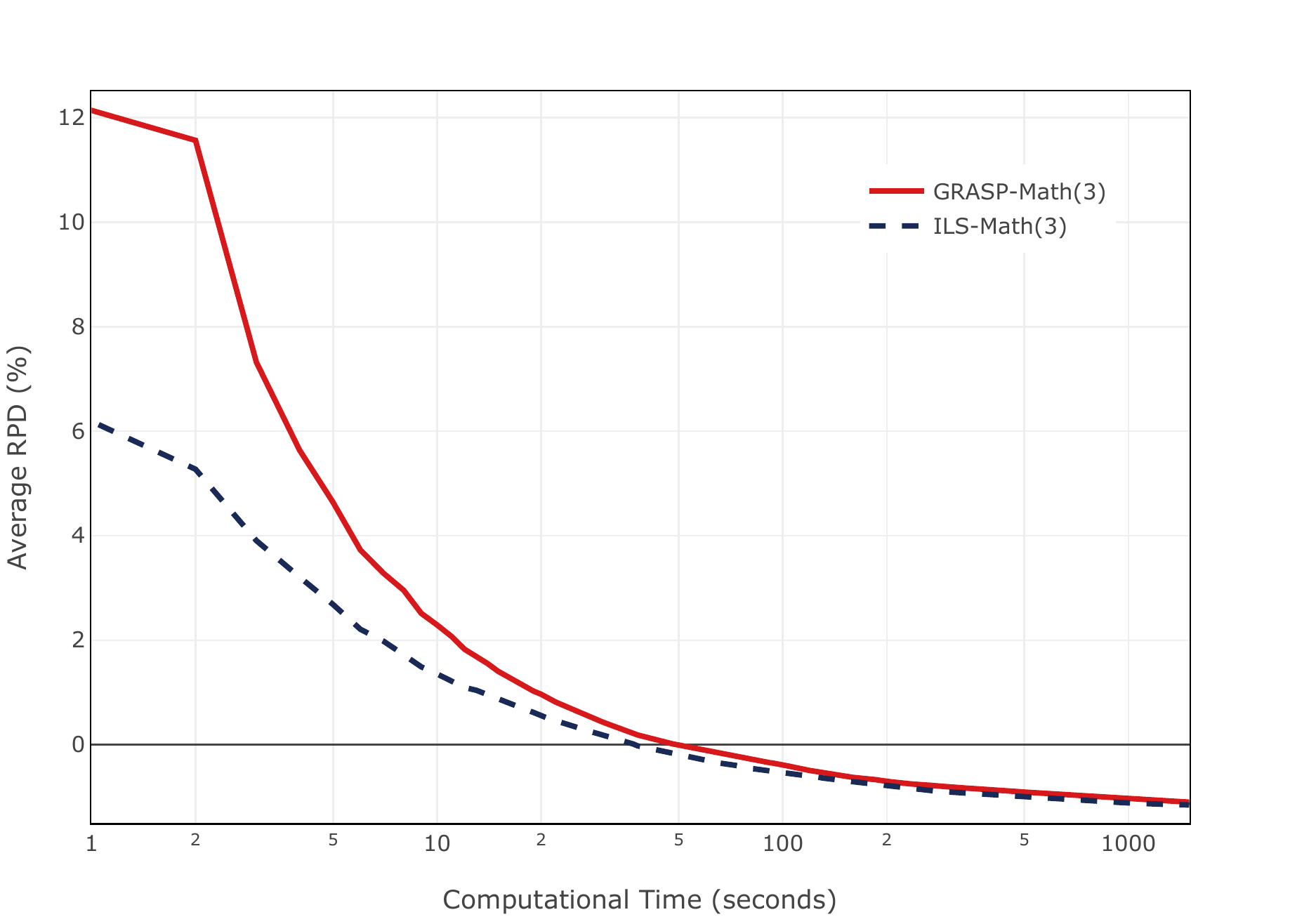}}
	\caption{Average RPD (vertical axis) evolution over the computational time (horizontal axis -- $\text{log}_{10}$ scale) for each matheuristic, running on the complete set of 72 PLSVSP instances. \label{fig:evolution}}
\end{figure}

As we evaluate solutions for a realistic process, it is important to analyze the final solutions of the proposed methods and their evolution over time, in case the company needs faster solutions. Based on this, we show in Figure~\ref{fig:evolution} the evolution of the solution for the two best approaches, \MHc\, and \MHf, illustrating the trade-off between the quality of the solution and the time spent for achieving it. For better visualization, the computational time axis is on a $\log_{10}$ scale. As the deviations and computational times are low for small-sized instances, the curve's shape is more defined by the medium and large-sized instances.

Note that both matheuristics surpass the BKS on average (that is, the current method used to solve the PLSVSP), crossing the zero line, in less than 1 minute (38 seconds for \MHc\,, and 49 seconds for \MHf). The \MHc\, matheuristic dominates the \MHf\, but reaching an equivalent average RPD after crossing the zero line. This analysis reinforces the advantages of using the matheuristics concept instead of running the MIP formulations, showing that the additional time-limit criteria can be included, if necessary, without significantly impacting the quality of the solution.

\section{Conclusions}
\label{sec:conclusion}

In this paper, we studied a ship scheduling problem related to offshore oil \& gas logistics, in which a company needs to schedule its fleet of Pipe Laying Support Vessels (PLSVs). These vessels are responsible for connecting sub-sea oil wells to the platforms, allowing production to begin. The company intends to anticipate the start of production in the most promising wells. The \textit{PLSV Scheduling Problem} (PLSVSP) can be seen as a variant of an identical parallel machine scheduling problem with batching and non-anticipatory family setup times to minimize the total weighted completion time. In this analogy, the vessels are the machines, the jobs are the wells, and each job includes several connection operations to be performed.

We introduced an ILS and a GRASP matheuristics, using two MIP-based neighborhood searches and a constructive heuristic to solve the PLSVSP. Two MIP formulations are considered, resulting in three variants of each method. The first is a batch formulation, which uses a WSPT dispatching rule to sequence operations within batches, while the second is a new formulation that considers sequencing within batches as a model decision. The results show that the matheuristics outperform the pure mathematical programming models in terms of computational time and solution quality. Among the matheuristics, a small advantage can be observed for two variants that combine the use of two batch formulations for the PLSVSP. The analysis shows that our new proposed formulation with sequencing variables helps the matheuristics to improve the quality of the solution. Moreover, new best solutions are provided for 40 of the 46 possible instances (without proven optimal solutions) on the PLSVSP instances benchmark.

The present study reinforces the importance of hybrid methods and their applicability in practical and theoretical contexts. The concept of splitting the problem to solve sub-problems is a relevant approach, especially when the size of the problem increases. It is worth mentioning that the approach is closely related to the nature of the PLSVSP in the studied company. In this environment, management guidelines for dealing with the problem change rapidly due to political and operational issues, requiring a quick adjustment in decision support tools. The use of mathematical formulation facilitates these adjustments, allowing the inclusion or exclusion of constraints without the need for substantial computational development in the algorithms. Nevertheless, regular metaheuristics based on local search procedures can be an interesting research direction due to the complexity of the problem, allowing other researchers working with similar or simplified versions of the PLSVSP to benefit from these studies. Moreover, the offshore environment where the PLSVs operate is surrounded by uncertainties caused by several aspects, such as climate changes (affecting ocean conditions), the complexity of operations, crew experience, and others, with significant impacts on the operating companies' planning. Thus, another attractive field for future research would be to understand the most impacted aspects of the problem and design stochastic optimization tools for the PLSVSP with uncertainties.

 \section*{Acknowledgements}

This study was financed in part by PUC-Rio, by the Coordenação de Aperfeiçoamento de Pessoal de Nível Superior – Brasil (CAPES) – Finance Code 001 and by the Conselho Nacional de Desenvolvimento Científico e Tecnológico (CNPq) under Grant numbers 403863/2016-3, 306802/2015-5, 425962/2016-4 and 313521/2017-4.

\appendix

\section{Parallel machine and PLSVSP definitions}
\label{app:mapping}

This section presents a mapping between the addressed identical parallel machine scheduling problem and its correspondence with the PLSVSP aspects (Table~\ref{tab:mapping}).

\definecolor{Gray}{gray}{0.9}

\begin{center}
	\begin{table}[htbp]
		\centering
		\footnotesize
		\caption{The mapping between the parallel machine scheduling problem definitions and its correspondences on the PLSVSP.}	
		\begin{tabularx}{\textwidth}{lXX}
			\toprule
			Name  & Machine scheduling definition & PLSVSP correspondence \\
			\midrule
			$\cO$     &  Operations  &  Pipeline connections \\
			$\cN$     &  Jobs  &  Wells\\
			$\cM$     &  Machines  &  PLSVs \\
			$\cF$     &  Families  & Groups of pipeline connections with similar loading process at the port \\
			$\cM_i$  &  Machine eligibility subset  &  Subset of vessels eligible to execute a pipeline connection \\
			$\cO_k$  &  Subset of operations a machine is eligible to execute  &  Subset of connections a vessel can carry \\
			$\cO_j$  &  Subset of operations composing a job &  Subset of connections required to enable a well to start producing \\
			$\cN_i$  &  Subset of jobs associated to an operation &  Subset of wells that depends on a connection to be able to produce \\
			$p_i$    &  Processing time of an operation &  Time taken to perform a pipeline connection \\
			$r_i$    &  Release date of an operation &  Arrival date of a pipeline at the port \\
			$l_i$    & Load occupancy of an operation &  Pipeline occupancy on the deck of the ship \\
			$f_i$  &  Family of an operation &  Group of a pipeline connection \\
			$r_k$  &  Release date of a machine &  Availability date of a vessel \\
			$q_k$  &  Capacity of a machine &  Vessel's deck capacity \\
			$w_j$  &  Weight of a job  &  Production potential of a well \\
			$C_i$  &  Completion time of an operation & Time when a connection ends\\
			$C_j$  &  Completion time of a job &  Time when a well is fully connected and able to start producing \\
			$s_f$  &  Setup times of a family &  Time spent at port loading pipelines for a specific group of operations \\
			$Batch$  &  Sequence of operations from the same family sharing the same setup time  &  PLSV voyage \\
			\bottomrule			
		\end{tabularx}%
		
		\label{tab:mapping}%
	\end{table}%
\end{center}

\section{Details of the example instance}
\label{app:example}

The details of the instance used to present a PLSV schedule example (Figure \ref{fig:scheduling-example}), are given below.
The details of the instance used to present a PLSV schedule example (Figure \ref{fig:scheduling-example}), are given below. The instance is composed by 15 operations ($|\cO|=15$), 5 jobs ($|\cN|=5$) and 4 machines ($|\cM|=4$). Table \ref{tab:opparameters} shows the processing time ($p_i$), the release date ($r_i$), the family ($f_i$), the load occupancy ($l_i$), the set of associated jobs ($\cN_i$), and the eligible set of machines ($\cM_i$) of each operation $i$. Table \ref{tab:jobparameters} presents the weight ($w_j$), and the subset of associated operations ($\cO_j$) of each job $j$. Table \ref{tab:machparameters} shows the release date ($r_k$), and the capacity ($q_k$) of each machine $k$, together with the subset of operations ($\cO_k$) that each machine $k$ is eligible to execute. Table \ref{tab:famparameters} presents the setup time ($s_f$) of each family $f$.

\begin{table}[htpb]
	\centering
	\caption{Operations aspects regarding the example instance.}
	\footnotesize
	\begin{tabular}{ccccccccccccc}
		\toprule
		\textbf{Operation ($i$)} &       & \textbf{$p_i$} &       & \textbf{$r_i$} &       & \textbf{$f_i$} &       & \textbf{$l_i$} &       & \textbf{$\cN_i$} &       & \textbf{$\cM_i$} \\
		\midrule
\textbf{1} &       & 23    &       & 5     &       & 1     &       & 90    &       & \{5\} &       & \{{1,4}\} \\
		\textbf{2} &       & 5     &       & 17    &       & 1     &       & 90    &       & \{5\} &       & \{1,3,4\} \\
		\textbf{3} &       & 25    &       & 16    &       & 3     &       & 70    &       & \{2,3,4\} &       & \{2,3,4\} \\
		\textbf{4} &       & 10    &       & 9     &       & 2     &       & 50    &       & \{1\} &       & \{1,2,3,4\} \\
		\textbf{5} &       & 4     &       & 18    &       & 1     &       & 60    &       & \{3\} &       & \{1,2,3,4\} \\
		\textbf{6} &       & 8     &       & 17    &       & 2     &       & 80    &       & \{1,3\} &       & \{2,4\} \\
		\textbf{7} &       & 29    &       & 15    &       & 3     &       & 60    &       & \{5\} &       & \{1,2,3,4\} \\
		\textbf{8} &       & 8     &       & 5     &       & 3     &       & 0     &       & \{2\} &       & \{2\} \\
		\textbf{9} &       & 17    &       & 9     &       & 3     &       & 40    &       & \{1,5\} &       & \{1,2,4\} \\
		\textbf{10} &       & 25    &       & 0     &       & 1     &       & 40    &       & \{2,5\} &       & \{1,2,3,4\} \\
		\textbf{11} &       & 4     &       & 3     &       & 1     &       & 30    &       & \{3\} &       & \{1,2,3,4\} \\
		\textbf{12} &       & 15    &       & 0     &       & 1     &       & 20    &       & \{1,3\} &       & \{1,3,4\} \\
		\textbf{13} &       & 7     &       & 16    &       & 2     &       & 60    &       & \{3\} &       & \{3,4\} \\
		\textbf{14} &       & 7     &       & 7     &       & 1     &       & 90    &       & \{4,5\} &       & \{1,4\} \\
		\textbf{15} &       & 18    &       & 15    &       & 3     &       & 20    &       & \{3\} &       & \{4\} \\
		\bottomrule
	\end{tabular}
	\label{tab:opparameters}%
\end{table}%

\begin{table}[htbp]
	\centering
	\footnotesize
	\caption{Jobs aspects regarding the example instance.}
	\begin{tabular}{ccccc}
		\toprule
		\textbf{Job ($j$)} &       & $w_j$  &       & $\cO_j$ \\
		\midrule
		\textbf{1} &       & 46    &       & \{6,9,11,12\} \\
		\textbf{2} &       & 40    &       & \{3,8,10\} \\
		\textbf{3} &       & 39    &       & \{3,4,5,6,12,13,15\} \\
		\textbf{4} &       & 13    &       & \{3,14\} \\
		\textbf{5} &       & 3     &       & \{1,2,7,9,10,14\} \\
		\bottomrule
	\end{tabular}
	\label{tab:jobparameters}%
\end{table}%

\begin{table}[htbp]
	\centering
	\footnotesize
	\caption{Machines aspects regarding the example instance.}
	\begin{tabular}{cccccc}
		\toprule
		\textbf{Machine ($k$)} &       & $r_k$  &       & $q_k$ & $\cO_k$\\
		\midrule
		\textbf{1} &       & 13    &       & 90 & \{1,2,4,5,7,9,10,11,12,14\}\\
		\textbf{2} &       & 0     &       & 80 & \{3,4,5,6,7,8,9,10,11\}\\
		\textbf{3} &       & 0     &       & 90 & \{2,3,4,5,7,10,11,12,13\}\\
		\textbf{4} &       & 1     &       & 90 & \{1,2,3,4,5,6,7,9,10,11,12,13,14,15\}\\
		\bottomrule
	\end{tabular}
	\label{tab:machparameters}%
\end{table}%

\begin{table}[htbp]
	\centering
	\caption{Family setup time duration regarding the example instance.}
	\begin{tabular}{ccc}
		\toprule
		\textbf{Family ($f$)} &       & $s_f$ \\
		\midrule
		\textbf{1} &       & 5 \\
		\textbf{2} &       & 7 \\
		\textbf{3} &       & 9 \\
		\bottomrule
	\end{tabular}
	\label{tab:famparameters}%
\end{table}%

\section{Matheuristics Comparative Analysis}
\label{app:mathcomp}

In this section, we compare the results for the matheuristic variations. In Table \ref{tab:ilsmath}, the results for the \texttt{ILS-Math} variations are shown in terms of the average RPD ($\overline{RPD}$) and average computational time ($\overline{time}$). The same results but for the \texttt{GRASP-Math} variations are presented in Table \ref{tab:graspmath}.

Note that results for the \texttt{ILS-Math} are quite similar for all variations. However, we can see an advantage for the \MHc\, approach when the size of the instances increases (instances with 25 and 50 operations). Regarding the complete set of instances, \MHc\, advantage is clearer with it dominating the other approaches in terms of the $\overline{RPD}$ with $-1.16\%$.

Similar behavior can be noted for the \texttt{GRASP-Math} variations. \MHf\, presented the best overall $\overline{RPD}$ of -1.119\% with a small difference in the average computational time for the best approach in this criterion (465 seconds against 456 seconds for the \MHd).

{\setlength{\tabcolsep}{6.7pt}
\begin{table}[htbp]
  \centering
  \footnotesize
    \caption{Average results for \texttt{ILS-Math} variations.}
    \begin{tabular}{ccccccccccc}
    \toprule
    \multirow{2}[4]{*}{$|\cO|$} & \multirow{2}[4]{*}{$|\cM|$} &       & \multicolumn{3}{c}{\MHa} & \multicolumn{2}{c}{\MHb} &       & \multicolumn{2}{c}{\MHc} \\
\cmidrule{4-5}\cmidrule{7-8}\cmidrule{10-11}          &       &       & $\overline{RPD}$   & $\overline{time}$   &       & $\overline{RPD}$   & $\overline{time}$   &       & $\overline{RPD}$   & $\overline{time}$ \\
\midrule
   \multirow{2}[1]{*}{15} & 4     &       & 0.61  & \textbf{14} &       & \textbf{0.33} & 18    &       & 0.42  & \textbf{14} \\
          & 8     &       & 0.11  & \textbf{46} &       & \textbf{0.03} & 49    &       & 0.08  & \textbf{46} \\
          &       &       &       &       &       &       &       &       &       &  \\
    \multirow{2}[0]{*}{25} & 4     &       & 0.06  & \textbf{123} &       & -0.13 & 143   &       & \textbf{-0.17} & \textbf{123} \\
          & 8     &       & -0.15 & 247   &       & -0.18 & 267   &       & \textbf{-0.20} & \textbf{244} \\
          &       &       &       &       &       &       &       &       &       &  \\
    \multirow{2}[1]{*}{50} & 4     &       & -2.78 & 1054  &       & -2.74 & \textbf{984} &       & \textbf{-3.02} & 1052 \\
          & 8     &       & -3.99 & \textbf{972} &       & -3.86 & 994   &       & \textbf{-4.09} & 990 \\
    \midrule
    \multicolumn{2}{c}{\multirow{2}[2]{*}{All instances}} &       &       &       &       &       &       &       &       &  \\
    \multicolumn{2}{c}{} &       & -1.02 & \textbf{409} &       & -1.09 & \textbf{409} &       & \textbf{-1.16} & 412 \\
    \bottomrule
    \end{tabular}%
    
  \label{tab:ilsmath}%
\end{table}%
}

{\setlength{\tabcolsep}{6pt}
\begin{table}[htbp]
  \centering
  \footnotesize
  \caption{Average results for \texttt{GRASP-Math} variations.}
    \begin{tabular}{ccrcccccccc}
    \toprule
    \multirow{2}[4]{*}{$|\cO|$} & \multirow{2}[4]{*}{$|\cM|$} &       & \multicolumn{3}{c}{\MHd} & \multicolumn{2}{c}{\MHe} &       & \multicolumn{2}{c}{\MHf} \\
\cmidrule{4-5}\cmidrule{7-8}\cmidrule{10-11}          &       &       & $\overline{RPD}$   & $\overline{time}$   &       & $\overline{RPD}$   & $\overline{time}$   &       & $\overline{RPD}$   & $\overline{time}$ \\
\midrule
       15    & 4     &       & 0.49  & \textbf{25} &       & \textbf{0.29} & 28    &       & 0.30  & 26 \\
    15    & 8     &       & 0.05  & \textbf{53} &       & 0.02  & 65    &       & \textbf{0.01} & 54 \\
          &       &       &       &       &       &       &       &       &       &  \\
    25    & 4     &       & -0.15 & \textbf{184} &       & \textbf{-0.38} & 217   &       & -0.31 & 186 \\
    25    & 8     &       & -0.20 & 298   &       & \textbf{-0.27} & 313   &       & \textbf{-0.27} & \textbf{297} \\
          &       &       &       &       &       &       &       &       &       &  \\
    50    & 4     &       & -2.34 & \textbf{1140} &       & -2.50 & 1194  &       & \textbf{-2.65} & 1159 \\
    50    & 8     &       & -3.74 & \textbf{1038} &       & \textbf{-3.86} & 1173  &       & -3.80 & 1066 \\
    \midrule
    \multicolumn{2}{c}{\multirow{2}[2]{*}{All instances}} &       &       &       &       &       &       &       &       &  \\
    \multicolumn{2}{c}{} &       & -0.98 & \textbf{456} &       & -1.117 & 498   &       & \textbf{-1.119} & 465 \\
    \bottomrule
    \end{tabular}%
    
  \label{tab:graspmath}%
\end{table}%
}


\section{Results by Instance}
\label{app:comp_results}

In this section, we present the total weighted completion time (objective function) found by each method for the benchmark of 72 PLSVSP instances with the respective computational times (Table \ref{tab:comp_results}). The first two columns indicate the name of the instance and its BKS from the literature (optimal solutions are indicated with an asterisk), defined by \cite{AbuMArrul2020}. Then, the upper bound (ub) and the computational time (time) for each formulation are presented, followed by the results of each matheuristic variation. The matheuristic results are presented in the following order: minimum total weighted completion time (min), average total weighted completion time (avg), and average computational time ($\overline{time}$), among the runs. We highlighted, in bold, the best solution for each instance. We suppress computational times that have reached the time limit (21,600 seconds) defined for the MIP formulations.

To save space, we shorten the instance names. For example, in the benchmark set, an instance named \textit{PLSV\_o15\_n5\_q3\_m4\_111\_1}, indicates a total of 15 operations (\textit{o15}), associated with 5 jobs (\textit{n5}), divided into 3 families (\textit{q3}), to be scheduled on 4 machines (\textit{m4}). Since the number of jobs depends on the number of operations (calculated as $|\cN| = \lfloor |\cO|/3 \rfloor$, following the instance generation procedure defined by \cite{AbuMArrul2020}), we delete it from the instance name. We also suppress the number of families and the last digit because they are equal to 3 and 1 in all instances.

\begin{landscape}

{\setlength{\tabcolsep}{1.4pt}
\renewcommand{\arraystretch}{1.15}
\tiny
\begin{longtable}{ccrrrrrrrrrrrrrrrrrrrrrrrrrrrrrrr}
\caption{Complete results by instance regarding the PLSVSP benchmark.} \label{tab:long} \\

\toprule
    \multirow{2}[4]{*}{Instance} &       & \multicolumn{1}{c}{\multirow{2}[4]{*}{BKS}} &       & \multicolumn{2}{c}{\texttt{Batch-WSPT}} &       & \multicolumn{2}{c}{\texttt{Batch-S}} &       & \multicolumn{3}{c}{\MHa} &       & \multicolumn{3}{c}{\MHb} &       & \multicolumn{3}{c}{\MHc} &       & \multicolumn{3}{c}{\MHd} &       & \multicolumn{3}{c}{\MHe} &       & \multicolumn{3}{c}{\MHf} \\
\cmidrule{5-6}\cmidrule{8-9}\cmidrule{11-13}\cmidrule{15-17}\cmidrule{19-21}\cmidrule{23-25}\cmidrule{27-29}\cmidrule{31-33}          &       &       &       & \multicolumn{1}{c}{ub} & \multicolumn{1}{c}{time} &       & \multicolumn{1}{c}{ub} & \multicolumn{1}{c}{time} &       & \multicolumn{1}{c}{min} & \multicolumn{1}{c}{avg} & \multicolumn{1}{c}{$\overline{time}$} &       & \multicolumn{1}{c}{min} & \multicolumn{1}{c}{avg} & \multicolumn{1}{c}{$\overline{time}$} &       & \multicolumn{1}{c}{min} & \multicolumn{1}{c}{avg} & \multicolumn{1}{c}{$\overline{time}$} &       & \multicolumn{1}{c}{min} & \multicolumn{1}{c}{avg} & \multicolumn{1}{c}{$\overline{time}$} &       & \multicolumn{1}{c}{min} & \multicolumn{1}{c}{avg} & \multicolumn{1}{c}{$\overline{time}$} &       & \multicolumn{1}{c}{min} & \multicolumn{1}{c}{avg} & \multicolumn{1}{c}{$\overline{time}$} \\
    \midrule
\endfirsthead

\multicolumn{15}{l}%
{{\bfseries \tablename\ \thetable{} -- continued from previous page}} \\
\toprule
    \multirow{2}[4]{*}{Instance} &       & \multicolumn{1}{c}{\multirow{2}[4]{*}{BKS}} &       & \multicolumn{2}{c}{\texttt{Batch-WSPT}} &       & \multicolumn{2}{c}{\texttt{Batch-S}} &       & \multicolumn{3}{c}{\MHa} &       & \multicolumn{3}{c}{\MHb} &       & \multicolumn{3}{c}{\MHc} &       & \multicolumn{3}{c}{\MHd} &       & \multicolumn{3}{c}{\MHe} &       & \multicolumn{3}{c}{\MHf} \\
\cmidrule{5-6}\cmidrule{8-9}\cmidrule{11-13}\cmidrule{15-17}\cmidrule{19-21}\cmidrule{23-25}\cmidrule{27-29}\cmidrule{31-33}          &       &       &       & \multicolumn{1}{c}{ub} & \multicolumn{1}{c}{time} &       & \multicolumn{1}{c}{ub} & \multicolumn{1}{c}{time} &       & \multicolumn{1}{c}{min} & \multicolumn{1}{c}{avg} & \multicolumn{1}{c}{$\overline{time}$} &       & \multicolumn{1}{c}{min} & \multicolumn{1}{c}{avg} & \multicolumn{1}{c}{$\overline{time}$} &       & \multicolumn{1}{c}{min} & \multicolumn{1}{c}{avg} & \multicolumn{1}{c}{$\overline{time}$} &       & \multicolumn{1}{c}{min} & \multicolumn{1}{c}{avg} & \multicolumn{1}{c}{$\overline{time}$} &       & \multicolumn{1}{c}{min} & \multicolumn{1}{c}{avg} & \multicolumn{1}{c}{$\overline{time}$} &       & \multicolumn{1}{c}{min} & \multicolumn{1}{c}{avg} & \multicolumn{1}{c}{$\overline{time}$} \\
    \midrule
\endhead
\hline \multicolumn{20}{l}{{Continued on next page}} \\
\endfoot

\bottomrule \\ \multicolumn{20}{l}{$\dagger$ Execution interrupted before 21,600 seconds (time limit) due to memory limit.}
\endlastfoot

    15-4-111 &       & 6903* &       & 6927  & 12417 &       & \textbf{6903} & 12121 &       & 6927  & 6930.8 & 2     &       & \textbf{6903} & 6903.0 & 4     &       & \textbf{6903} & 6909.2 & 2     &       & \textbf{6903} & 6915.0 & 5     &       & \textbf{6903} & 6903.0 & 4     &       & \textbf{6903} & 6903.0 & 5 \\
    15-4-112 &       & 7717*  &       & 7842  & -     &       & \textbf{7717} & -     &       & 7760  & 7763.6 & 21    &       & \textbf{7717} & 7755.3 & 26    &       & 7760  & 7763.6 & 21    &       & \textbf{7717} & 7742.6 & 42    &       & \textbf{7717} & 7736.2 & 46    &       & \textbf{7717} & 7742.6 & 43 \\
    15-4-121 &       & 7498* &       & 7589  & -     &       & \textbf{7498} & -     &       & 7589  & 7609.6 & 5     &       & \textbf{7498} & 7539.6 & 8     &       & \textbf{7498} & 7532.4 & 5     &       & 7589  & 7595.1 & 7     &       & \textbf{7498} & 7500.7 & 10    &       & \textbf{7498} & 7498.0 & 8 \\
    15-4-122 &       & 8608* &       & \textbf{8608} & -     &       & \textbf{8608} & -     &       & \textbf{8608} & 8608.0 & 10    &       & \textbf{8608} & 8612.0 & 15    &       & \textbf{8608} & 8608.0 & 10    &       & \textbf{8608} & 8608.0 & 18    &       & \textbf{8608} & 8608.0 & 26    &       & \textbf{8608} & 8608.0 & 18 \\
    15-4-131 &       & 10964* &       & \textbf{10964} & 16546 &       & \textbf{10964} & 20576 &       & \textbf{10964} & 10996.0 & 2     &       & \textbf{10964} & 10983.0 & 2     &       & \textbf{10964} & 10996.0 & 2     &       & \textbf{10964} & 11222.5 & 0     &       & \textbf{10964} & 11292.5 & 0     &       & \textbf{10964} & 11222.5 & 0 \\
    15-4-132 &       & 10685* &       & \textbf{10685} & -     &       & \textbf{10685} & -     &       & \textbf{10685} & 10685.0 & 5     &       & \textbf{10685} & 10685.0 & 6     &       & \textbf{10685} & 10685.0 & 5     &       & \textbf{10685} & 10687.1 & 10    &       & \textbf{10685} & 10687.1 & 11    &       & \textbf{10685} & 10687.1 & 10 \\
    15-4-211 &       & 8822* &       & \textbf{8822} & -     &       & \textbf{8822} & -     &       & \textbf{8822} & 8872.0 & 12    &       & \textbf{8822} & 8859.5 & 11    &       & \textbf{8822} & 8872.0 & 12    &       & \textbf{8822} & 8822.0 & 24    &       & \textbf{8822} & 8834.5 & 19    &       & \textbf{8822} & 8822.0 & 24 \\
    15-4-212 &       & 6579* &       & 6681  & -     &       & \textbf{6579} & -     &       & 6681  & 6684.6 & 49    &       & \textbf{6579} & 6680.8 & 61    &       & \textbf{6579} & 6628.2 & 51    &       & \textbf{6579} & 6611.0 & 104   &       & \textbf{6579} & 6579.0 & 89    &       & \textbf{6579} & 6579.0 & 109 \\
    15-4-221 &       & 5929* &       & \textbf{5929} & -     &       & \textbf{5929} & -     &       & \textbf{5929} & 5929.0 & 11    &       & \textbf{5929} & 5929.0 & 14    &       & \textbf{5929} & 5929.0 & 11    &       & \textbf{5929} & 5929.0 & 18    &       & \textbf{5929} & 5929.0 & 25    &       & \textbf{5929} & 5929.0 & 19 \\
    15-4-222 &       & 14545* &       & 14712 & -     &       & \textbf{14545} & -     &       & 14712 & 14754.6 & 18    &       & \textbf{14545} & 14621.9 & 21    &       & 14712 & 14754.6 & 18    &       & \textbf{14545} & 14628.8 & 27    &       & \textbf{14545} & 14545.0 & 40    &       & \textbf{14545} & 14628.8 & 27 \\
    15-4-231 &       & 10755* &       & \textbf{10755} & -     &       & \textbf{10755} & -     &       & \textbf{10755} & 10859.3 & 15    &       & \textbf{10755} & 10781.5 & 25    &       & \textbf{10755} & 10845.7 & 16    &       & \textbf{10755} & 10817.6 & 26    &       & \textbf{10755} & 10764.1 & 36    &       & \textbf{10755} & 10790.4 & 26 \\
    15-4-232 &       & 15747* &       & \textbf{15747} & -     &       & \textbf{15747} & -     &       & \textbf{15747} & 15747.0 & 13    &       & \textbf{15747} & 15747.0 & 17    &       & \textbf{15747} & 15747.0 & 13    &       & \textbf{15747} & 15758.2 & 22    &       & \textbf{15747} & 15747.0 & 28    &       & \textbf{15747} & 15747.0 & 23 \\
    \midrule
    15-8-111 &       & 4306* &       & \textbf{4306} & -     &       & \textbf{4306} & -     &       & \textbf{4306} & 4306.0 & 32    &       & \textbf{4306} & 4306.0 & 42    &       & \textbf{4306} & 4306.0 & 32    &       & \textbf{4306} & 4306.0 & 44    &       & \textbf{4306} & 4306.0 & 56    &       & \textbf{4306} & 4306.0 & 45 \\
    15-8-112 &       & 5453* &       & \textbf{5453} & -     &       & \textbf{5453} & -     &       & \textbf{5453} & 5453.0 & 73    &       & \textbf{5453} & 5453.0 & 97    &       & \textbf{5453} & 5453.0 & 73    &       & \textbf{5453} & 5453.0 & 92    &       & \textbf{5453} & 5453.0 & 117   &       & \textbf{5453} & 5453.0 & 93 \\
    15-8-121 &       & 10790* &       & \textbf{10790} & 0     &       & \textbf{10790} & 1     &       & \textbf{10790} & 10790.0 & 1     &       & \textbf{10790} & 10790.0 & 2     &       & \textbf{10790} & 10790.0 & 1     &       & \textbf{10790} & 10790.0 & 1     &       & \textbf{10790} & 10790.0 & 2     &       & \textbf{10790} & 10790.0 & 1 \\
    15-8-122 &       & 8369* &       & \textbf{8369} & 2     &       & \textbf{8369} & 2     &       & \textbf{8369} & 8369.0 & 9     &       & \textbf{8369} & 8369.0 & 10    &       & \textbf{8369} & 8369.0 & 9     &       & \textbf{8369} & 8369.0 & 19    &       & \textbf{8369} & 8369.0 & 17    &       & \textbf{8369} & 8369.0 & 19 \\
    15-8-131 &       & 4339* &       & \textbf{4339} & -     &       & \textbf{4339} & -     &       & \textbf{4339} & 4339.0 & 10    &       & \textbf{4339} & 4339.0 & 13    &       & \textbf{4339} & 4339.0 & 10    &       & \textbf{4339} & 4339.0 & 14    &       & \textbf{4339} & 4339.0 & 15    &       & \textbf{4339} & 4339.0 & 14 \\
    15-8-132 &       & 7371* &       & \textbf{7371} & -     &       & \textbf{7371} & -     &       & \textbf{7371} & 7386.3 & 125   &       & \textbf{7371} & 7379.8 & 123   &       & \textbf{7371} & 7386.3 & 125   &       & \textbf{7371} & 7382.3 & 132   &       & \textbf{7371} & 7383.1 & 136   &       & \textbf{7371} & 7382.9 & 132 \\
    15-8-211 &       & 3189* &       & 3203  & -     &       & \textbf{3189} & -     &       & 3203  & 3203.3 & 72    &       & \textbf{3189} & 3191.8 & 69    &       & \textbf{3189} & 3193.2 & 74    &       & 3203  & 3203.0 & 74    &       & \textbf{3189} & 3189.0 & 109   &       & \textbf{3189} & 3189.0 & 78 \\
    15-8-212 &       & 4256* &       & \textbf{4256} & -     &       & 4263  & -     &       & \textbf{4256} & 4256.7 & 112   &       & \textbf{4256} & 4258.5 & 97    &       & \textbf{4256} & 4256.7 & 112   &       & \textbf{4256} & 4256.0 & 111   &       & \textbf{4256} & 4257.4 & 133   &       & \textbf{4256} & 4256.0 & 111 \\
    15-8-221 &       & 5519* &       & \textbf{5519} & -     &       & \textbf{5519} & -     &       & \textbf{5519} & 5519.0 & 26    &       & \textbf{5519} & 5522.0 & 30    &       & \textbf{5519} & 5519.0 & 26    &       & \textbf{5519} & 5519.0 & 35    &       & \textbf{5519} & 5519.0 & 44    &       & \textbf{5519} & 5519.0 & 36 \\
    15-8-222 &       & 10461* &       & \textbf{10461} & -     &       & \textbf{10461} & -     &       & \textbf{10461} & 10530.3 & 7     &       & \textbf{10461} & 10461.0 & 8     &       & \textbf{10461} & 10530.3 & 7     &       & \textbf{10461} & 10461.0 & 9     &       & \textbf{10461} & 10461.0 & 11    &       & \textbf{10461} & 10461.0 & 9 \\
    15-8-231 &       & 8002* &       & \textbf{8002} & -     &       & \textbf{8002} & -     &       & \textbf{8002} & 8002.0 & 36    &       & \textbf{8002} & 8005.0 & 35    &       & \textbf{8002} & 8002.0 & 36    &       & \textbf{8002} & 8002.0 & 55    &       & \textbf{8002} & 8003.5 & 68    &       & \textbf{8002} & 8002.0 & 55 \\
    15-8-232 &       & 5127* &       & \textbf{5127} & -     &       & \textbf{5127} & -     &       & \textbf{5127} & 5127.0 & 50    &       & \textbf{5127} & 5127.0 & 60    &       & \textbf{5127} & 5127.0 & 51    &       & \textbf{5127} & 5127.0 & 54    &       & \textbf{5127} & 5127.0 & 69    &       & \textbf{5127} & 5127.0 & 54 \\
    \midrule
    25-4-111 &       & 9151  &       & 9192  & -     &       & \textbf{9151} & -     &       & 9192  & 9217.3 & 118   &       & 9166  & 9237.1 & 102   &       & \textbf{9151} & 9176.1 & 119   &       & \textbf{9151} & 9174.3 & 150   &       & \textbf{9151} & 9157.3 & 161   &       & \textbf{9151} & 9152.2 & 150 \\
    25-4-112 &       & 18776 &       & 18896 & -     &       & 19012 & -     &       & 18752 & 18795.2 & 151   &       & 18734 & 18817.7 & 174   &       & \textbf{18678} & 18751.8 & 153   &       & 18752 & 18790.9 & 192   &       & \textbf{18678} & 18707.4 & 222   &       & \textbf{18678} & 18743.5 & 195 \\
    25-4-121 &       & 23134 &       & 23364 & -     &       & 23365 & -     &       & 22927 & 23080.1 & 155   &       & \textbf{22865} & 23011.0 & 253   &       & \textbf{22865} & 23050.7 & 161   &       & \textbf{22865} & 23039.0 & 263   &       & \textbf{22865} & 23026.9 & 331   &       & \textbf{22865} & 22986.4 & 294 \\
    25-4-122 &       & 12415 &       & 12417 & -     &       & \textbf{12415} & -     &       & 12417 & 12417.6 & 51    &       & \textbf{12415} & 12415.6 & 52    &       & \textbf{12415} & 12416.8 & 51    &       & \textbf{12415} & 12429.9 & 77    &       & \textbf{12415} & 12415.0 & 90    &       & \textbf{12415} & 12422.8 & 78 \\
    25-4-131 &       & 32800 &       & 33092 & -     &       & 32966 & -     &       & 32854 & 32862.3 & 176   &       & \textbf{32800} & 32859.9 & 215   &       & 32854 & 32862.3 & 176   &       & \textbf{32800} & 32848.7 & 314   &       & 32820 & 32863.4 & 224   &       & \textbf{32800} & 32860.9 & 295 \\
    25-4-132 &       & 27555* &       & 28036 & -     &       & \textbf{27555} & -     &       & 27651 & 27663.2 & 93    &       & \textbf{27555} & 27608.2 & 109   &       & \textbf{27555} & 27555.0 & 88    &       & \textbf{27555} & 27633.0 & 188   &       & \textbf{27555} & 27562.4 & 234   &       & \textbf{27555} & 27555.0 & 187 \\
    25-4-211 &       & 30098 &       & 30534 & -     &       & 30301 & -     &       & 30002 & 30103.6 & 132   &       & \textbf{29679} & 29924.1 & 184   &       & 30002 & 30099.3 & 133   &       & 29759 & 29794.8 & 180   &       & \textbf{29679} & 29710.2 & 220   &       & \textbf{29679} & 29766.2 & 181 \\
    25-4-212 &       & 20012 &       & 19669 & -     &       & 20667 & -     &       & 19669 & 19700.3 & 129   &       & \textbf{19595} & 19639.4 & 146   &       & \textbf{19595} & 19638.8 & 131   &       & 19628 & 19719.2 & 188   &       & \textbf{19595} & 19607.4 & 310   &       & \textbf{19595} & 19694.9 & 191 \\
    25-4-221 &       & 19944 &       & 20399 & -     &       & 20637 & -     &       & 19944 & 19970.3 & 120   &       & \textbf{19833} & 19955.7 & 148   &       & \textbf{19833} & 19929.8 & 124   &       & 19944 & 20018.9 & 155   &       & \textbf{19833} & 19945.3 & 225   &       & \textbf{19833} & 19987.3 & 153 \\
    25-4-222 &       & 27274 &       & 27220 & -     &       & 27268 & -     &       & 27307 & 27360.2 & 60    &       & \textbf{27148} & 27267.7 & 107   &       & 27185 & 27262.1 & 61    &       & 27155 & 27219.0 & 80    &       & \textbf{27148} & 27165.9 & 162   &       & \textbf{27148} & 27184.7 & 77 \\
    25-4-231 &       & 29552 &       & 29695 & -     &       & \textbf{29552} & -     &       & 29695 & 29698.2 & 91    &       & \textbf{29552} & 29580.6 & 81    &       & \textbf{29552} & 29552.0 & 89    &       & \textbf{29552} & 29638.1 & 174   &       & \textbf{29552} & 29552.0 & 155   &       & \textbf{29552} & 29552.0 & 186 \\
    25-4-232 &       & 29502 &       & 29894 & -     &       & \textbf{29383} & -     &       & 29540 & 29552.0 & 201   &       & \textbf{29383} & 29431.2 & 144   &       & \textbf{29383} & 29472.7 & 186   &       & \textbf{29383} & 29407.1 & 251   &       & \textbf{29383} & 29408.8 & 273   &       & \textbf{29383} & 29400.8 & 245 \\
    \midrule
    25-8-111 &       & 11558 &       & 11634 & -     &       & 11479 & -     &       & \textbf{11387} & 11464.1 & 204   &       & \textbf{11387} & 11447.5 & 255   &       & \textbf{11387} & 11464.1 & 205   &       & 11427 & 11465.1 & 361   &       & \textbf{11387} & 11444.4 & 334   &       & 11427 & 11453.0 & 363 \\
    25-8-112 &       & 16150 &       & 16191 & -     &       & 16295 & -     &       & 16165 & 16267.8 & 309   &       & 16056 & 16203.4 & 273   &       & 16061 & 16180.9 & 300   &       & \textbf{16053} & 16123.8 & 373   &       & \textbf{16053} & 16095.7 & 389   &       & \textbf{16053} & 16089.4 & 354 \\
    25-8-121 &       & 10478 &       & \textbf{10478} & -     &       & 10571 & -     &       & \textbf{10478} & 10537.1 & 282   &       & \textbf{10478} & 10539.0 & 278   &       & \textbf{10478} & 10548.8 & 254   &       & \textbf{10478} & 10544.9 & 287   &       & \textbf{10478} & 10541.9 & 294   &       & \textbf{10478} & 10548.5 & 309 \\
    25-8-122 &       & 19723 &       & 19802 & -     &       & 19778 & -     &       & \textbf{19647} & 19658.0 & 154   &       & \textbf{19647} & 19659.4 & 187   &       & \textbf{19647} & 19658.0 & 150   &       & \textbf{19647} & 19661.5 & 177   &       & \textbf{19647} & 19668.6 & 206   &       & \textbf{19647} & 19659.5 & 177 \\
    25-8-131 &       & 9700* &       & 9716  & -     &       & \textbf{9700} & -     &       & \textbf{9700} & 9700.0 & 145   &       & \textbf{9700} & 9702.2 & 151   &       & \textbf{9700} & 9700.0 & 145   &       & \textbf{9700} & 9700.9 & 183   &       & \textbf{9700} & 9700.6 & 184   &       & \textbf{9700} & 9700.0 & 185 \\
    25-8-132 &       & 17947 &       & \textbf{17911} & -     &       & 17989 & -     &       & \textbf{17911} & 17931.0 & 198   &       & \textbf{17911} & 17941.6 & 223   &       & \textbf{17911} & 17934.8 & 194   &       & \textbf{17911} & 17964.6 & 214   &       & \textbf{17911} & 17946.4 & 265   &       & \textbf{17911} & 17956.2 & 225 \\
    25-8-211 &       & 8369  &       & 8367  & -     &       & 8429  & -     &       & \textbf{8261} & 8315.2 & 307   &       & \textbf{8261} & 8316.1 & 383   &       & \textbf{8261} & 8325.6 & 384   &       & 8308  & 8334.3 & 405   &       & \textbf{8261} & 8327.1 & 391   &       & 8294  & 8328.9 & 372 \\
    25-8-212 &       & 13518 &       & 13675 & -     &       & 13963 & -     &       & 13397 & 13539.1 & 298   &       & 13388 & 13505.4 & 378   &       & 13397 & 13530.6 & 275   &       & \textbf{13337} & 13472.4 & 432   &       & \textbf{13337} & 13454.3 & 417   &       & \textbf{13337} & 13454.1 & 413 \\
    25-8-221 &       & 10993 &       & 11026 & -     &       & 11004 & -     &       & 11019 & 11029.4 & 150   &       & \textbf{10993} & 11036.0 & 151   &       & \textbf{10993} & 11012.9 & 150   &       & 11004 & 11025.0 & 144   &       & \textbf{10993} & 11000.8 & 189   &       & \textbf{10993} & 11006.7 & 148 \\
    25-8-222 &       & 14140 &       & 13998 & -     &       & 14001 & -     &       & \textbf{13819} & 13937.4 & 345   &       & \textbf{13819} & 13973.6 & 360   &       & 13845 & 13934.4 & 327   &       & 13902 & 13968.9 & 314   &       & 13845 & 13947.0 & 422   &       & 13844 & 13939.0 & 365 \\
    25-8-231 &       & 9251  &       & 9312  & -     &       & 9282  & -     &       & \textbf{9230} & 9232.5 & 241   &       & \textbf{9230} & 9237.6 & 229   &       & \textbf{9230} & 9231.4 & 239   &       & \textbf{9230} & 9233.9 & 348   &       & \textbf{9230} & 9241.8 & 265   &       & \textbf{9230} & 9238.2 & 338 \\
    25-8-232 &       & 16929 &       & 17005 & -     &       & 17324 & -     &       & 16867 & 16913.5 & 325   &       & \textbf{16863} & 16916.8 & 337   &       & \textbf{16863} & 16906.2 & 304   &       & 16889 & 16930.2 & 333   &       & 16880 & 16937.9 & 395   &       & 16879 & 16932.3 & 316 \\
    \midrule
    50-4-111 &       & 64068 &       & 68673 & -     &       & 70886 & -     &       & 61288 & 61935.1 & 1381  &       & 61237 & 61943.4 & 1238  &       & 61179 & 61697.1 & 1159  &       & 61346 & 62766.2 & 1311  &       & 61577 & 62328.9 & 1541  &       & \textbf{61023} & 62129.3 & 1242 \\
    50-4-112 &       & 102246 &       & 102201 & -     &       & 102283 & -     &       & 98283 & 98719.3 & 1165  &       & 98116 & 98668.3 & 1127  &       & \textbf{97260} & 98259.4 & 1010  &       & 98093 & 98912.2 & 1511  &       & 98332 & 99184.0 & 1082  &       & 98582 & 99068.3 & 1245 \\
    50-4-121 &       & 65231 &       & 65185 & -     &       & 65367 & -     &       & 63447 & 63728.6 & 1075  &       & 63578 & 63947.6 & 820   &       & \textbf{63395} & 63648.4 & 1038  &       & 63736 & 64064.0 & 1170  &       & 63693 & 64103.4 & 1018  &       & 63752 & 64133.9 & 974 \\
    50-4-122 &       & 109078 &       & 110531 & -     &       & 108577 & -     &       & 106608 & 107269.2 & 1067  &       & 106194 & 107162.2 & 895   &       & \textbf{106009} & 106811.5 & 1418  &       & 106776 & 107705.8 & 1153  &       & 106388 & 107780.3 & 1044  &       & 106491 & 107224.0 & 1229 \\
    50-4-131 &       & 92825 &       & 92011 & -     &       & 92386 & -     &       & 91497 & 91748.6 & 1166  &       & \textbf{91341} & 91698.6 & 1075  &       & 91481 & 91751.8 & 831   &       & 91532 & 92005.7 & 1047  &       & 91566 & 91803.9 & 1034  &       & 91398 & 91675.4 & 1054 \\
    50-4-132 &       & 99080 &       & 99690 & -     &       & 99599 & -     &       & \textbf{97906} & 98372.1 & 818   &       & 98343 & 98593.0 & 803   &       & \textbf{97906} & 98448.1 & 601   &       & 98208 & 98858.8 & 852   &       & 98382 & 98782.8 & 1095  &       & 98208 & 98760.3 & 1109 \\
    50-4-211 &       & 84383 &       & 82439 & -     &       & 81604 & -     &       & 80032 & 80469.9 & 1061  &       & 80254 & 80686.1 & 1160  &       & 79436 & 80281.2 & 1110  &       & 79792 & 80779.7 & 1193  &       & 79611 & 80506.8 & 1390  &       & \textbf{78988} & 80475.2 & 1333 \\
    50-4-212 &       & 87761 &       & 83410 & -     &       & 86607 & -     &       & 81350 & 81920.5 & 805   &       & 81448 & 82028.5 & 1110  &       & 81145 & 81713.8 & 934   &       & 81881 & 82204.4 & 1210  &       & 81678 & 81995.5 & 1453  &       & \textbf{80798} & 81669.1 & 1308 \\
    50-4-221 &       & 106430 &       & 105371 & -     &       & 104586 & -     &       & 103930 & 104350.2 & 1011  &       & 103516 & 104005.3 & 1086  &       & \textbf{103345} & 104019.3 & 1125  &       & 103732 & 104368.3 & 1449  &       & 103496 & 104426.2 & 1373  &       & 103538 & 104292.5 & 1232 \\
    50-4-222 &       & 105136 &       & 104714 & 20784$^{\dagger}$ &       & 107833 & -     &       & 101723 & 102670.7 & 1177  &       & 101794 & 102557.0 & 768   &       & \textbf{100921} & 102190.6 & 1215  &       & 102113 & 102938.4 & 959   &       & 101563 & 102465.9 & 914   &       & 101225 & 102343.6 & 1234 \\
    50-4-231 &       & 89865 &       & 85913 & -     &       & 87968 & -     &       & 85588 & 85814.1 & 890   &       & 85401 & 85746.2 & 912   &       & \textbf{85290} & 85701.7 & 1062  &       & 85456 & 85827.3 & 797   &       & 85475 & 85839.8 & 1401  &       & 85393 & 85923.0 & 892 \\
    50-4-232 &       & 131034 &       & 133193 & 2948$^{\dagger}$  &       & 134424 & -     &       & 129575 & 130260.5 & 1032  &       & 129334 & 130449.9 & 808   &       & \textbf{128913} & 129813.0 & 1125  &       & 130160 & 131671.3 & 1028  &       & 129823 & 131231.3 & 979   &       & 129649 & 130932.7 & 1056 \\
    \midrule
    50-8-111 &       & 34381 &       & 33354 & -     &       & 35691 & -     &       & \textbf{33248} & 33601.0 & 989   &       & 33402 & 33669.0 & 1029  &       & 33269 & 33514.7 & 997   &       & 33672 & 33804.2 & 933   &       & 33590 & 33751.7 & 901   &       & 33574 & 33800.8 & 882 \\
    50-8-112 &       & 50521 &       & 47910 & -     &       & 48461 & -     &       & 45320 & 45815.3 & 1074  &       & 45596 & 46078.9 & 899   &       & \textbf{44980} & 45704.9 & 1177  &       & 45906 & 46149.1 & 1191  &       & 45965 & 46225.8 & 1521  &       & 45609 & 46067.1 & 1676 \\
    50-8-121 &       & 31944 &       & 30906 & -     &       & 32771 & 6026$^{\dagger}$  &       & 30393 & 30564.7 & 786   &       & \textbf{30276} & 30493.0 & 977   &       & 30349 & 30593.1 & 891   &       & 30332 & 30618.4 & 832   &       & 30321 & 30541.1 & 849   &       & 30336 & 30591.2 & 870 \\
    50-8-122 &       & 42961 &       & 41870 & -     &       & 43002 & -     &       & 41111 & 41419.7 & 890   &       & 41101 & 41496.5 & 875   &       & 41147 & 41379.2 & 1196  &       & 41178 & 41389.9 & 1077  &       & \textbf{41080} & 41339.6 & 1095  &       & 41204 & 41360.6 & 1173 \\
    50-8-131 &       & 60738 &       & 59588 & -     &       & 60242 & -     &       & \textbf{59018} & 59450.0 & 673   &       & 59377 & 59538.9 & 857   &       & 59111 & 59483.5 & 765   &       & 59276 & 59684.5 & 968   &       & 59410 & 59651.2 & 951   &       & 59335 & 59629.8 & 1096 \\
    50-8-132 &       & 65103 &       & 64953 & 17470$^{\dagger}$ &       & 68600 & 4164$^{\dagger}$  &       & 62657 & 63114.0 & 997   &       & 62797 & 63071.0 & 942   &       & 62486 & 62828.3 & 805   &       & \textbf{62387} & 63022.3 & 1118  &       & 62665 & 63116.7 & 988   &       & 62758 & 63019.4 & 1118 \\
    50-8-211 &       & 48355 &       & 48884 & -     &       & 50349 & -     &       & \textbf{46323} & 46927.1 & 1207  &       & 46664 & 46987.4 & 1258  &       & 46489 & 46861.7 & 1102  &       & 46589 & 47007.7 & 1400  &       & 46586 & 46971.6 & 1390  &       & 46510 & 47026.8 & 1178 \\
    50-8-212 &       & 60291 &       & 58338 & 6810$^{\dagger}$  &       & 60588 & 5296$^{\dagger}$  &       & 55215 & 55913.3 & 1258  &       & 55512 & 56081.6 & 1268  &       & 55375 & 55858.9 & 1417  &       & \textbf{54994} & 55895.1 & 1393  &       & 55372 & 55711.6 & 1586  &       & 55528 & 55873.1 & 1314 \\
    50-8-221 &       & 35471 &       & 36247 & 5467$^{\dagger}$  &       & 36632 & 5656$^{\dagger}$  &       & 33847 & 34086.6 & 1164  &       & 33894 & 34044.6 & 1138  &       & 33847 & 34081.7 & 883   &       & 33892 & 34208.0 & 897   &       & 33836 & 34141.8 & 1258  &       & \textbf{33711} & 34120.3 & 1207 \\
    50-8-222 &       & 56600 &       & 54512 & -     &       & 57001 & 13619$^{\dagger}$ &       & 52738 & 53360.1 & 928   &       & 53051 & 53521.2 & 960   &       & 52740 & 53395.5 & 1122  &       & 53281 & 53579.8 & 948   &       & \textbf{52603} & 53360.3 & 1509  &       & 53094 & 53642.8 & 804 \\
    50-8-231 &       & 54080 &       & 53303 & -     &       & 53941 & -     &       & 52805 & 53019.7 & 844   &       & 52756 & 53049.8 & 1004  &       & \textbf{52699} & 52913.3 & 831   &       & 52896 & 53115.6 & 1041  &       & 52791 & 53058.1 & 1311  &       & 52940 & 53091.0 & 843 \\
    50-8-232 &       & 62858 &       & 62385 & 6017$^{\dagger}$  &       & 67546 & 2451$^{\dagger}$  &       & \textbf{61568} & 61976.8 & 856   &       & 61852 & 62099.2 & 724   &       & 61606 & 62009.1 & 701   &       & 61865 & 62197.0 & 654   &       & 61713 & 62091.8 & 722   &       & 61598 & 62105.5 & 630 \\

    \label{tab:comp_results}
\end{longtable}
}
\end{landscape}


\end{document}